\documentclass[12pt,reqno]{amsart}
\usepackage[margin=1in]{geometry}
\usepackage{amsmath,amssymb,amsthm}
\usepackage{mathrsfs}
\usepackage{xcolor}
\usepackage{hyperref}
\allowdisplaybreaks
\usepackage{quotes}

\numberwithin{equation}{section}
\newtheorem{Thm}{Theorem}[section]
\newtheorem{Cor}{Corollary}[section]
\newtheorem{Def}{Definition}[section]
\newtheorem{Lem}{Lemma}[section]
\newtheorem{Pro}{Proposition}[section]

\newtheorem{Hyp}{Hypothesis}[section]

\def\mR{\mathbb{R}}
 
\def\mN{\mathbb{N}}

\def\mC{\mathbb{C}}

\def\mZ{\mathbb{Z}}

\def\mcF{\mathbb {\mathcal F}}

\def\mcA{\mathcal{A}}

\def\macP{\mathcal{P}}
\def\macQ{\mathcal{Q}}

\def\m{{\bf m}}
\def\n{{\bf n}}

\def\T{{\bf T}}

\def\i0t{\int_0^t}

\let\originalleft\left
\let\originalright\right
\renewcommand{\left}{\mathopen{}\mathclose\bgroup\originalleft}
\renewcommand{\right}{\aftergroup\egroup\originalright}

\begin{document}
	
\newcommand{\Addresses}{{
		\bigskip
		\footnote{
	\noindent  \textsuperscript{1} U. S. Air Force Research Laboratory, Wright Patterson Air Force Base, Ohio 45433, U. S. A.
	
	\par\nopagebreak
	\noindent  \textsuperscript{2} NRC-Senior Research Fellow, National Academies of Science, Engineering and Medicine,  U. S. Air Force Research Laboratory, Wright Patterson Air Force Base, Ohio 45433, U. S. A.		
		\par\nopagebreak \noindent
		\textit{e-mail:} \texttt{provostsritharan@gmail.com} 	$^*$Corresponding author.
	\par\nopagebreak
\noindent  \textit{e-mail:} \texttt{saba.mudaliar@us.af.mil}

}}}

\title[Nonlinear Filtering of Spin Systems]{Nonlinear Filtering of Classical and Quantum Spin Systems \Addresses}
\author[ S. S. Sritharan and Saba Mudaliar]
{  Sivaguru S. Sritharan\textsuperscript{1,2*} and Saba Mudaliar\textsuperscript{1}}

\maketitle
{\hspace{1.5in}\footnotesize\dedicatory{Dedicated in honor of Professor K. R. Parthasarathy}}

\begin{abstract}
In this paper we consider classical and quantum spin systems on discrete lattices and in Euclidean spaces, modeled by infinite dimensional stochastic diffusions in Hilbert spaces. Existence and uniqueness of various notions of solutions, existence and uniqueness of invariant measures as well as exponential convergence to equilibrium are known for these models. We formulate nonlinear filtering problem for these classes of models, derive nonlinear filtering equations of Fujisaki-Kallianpur-Kunita and Zakai tye, and prove existence and uniqueness of measure-valued solutions to these filtering equations. We then establish the Feller property and Markov property of the semigroups associated with the filtering equations and also prove existence and uniqueness of invariant measures. Evolution of error covariance equation for the nonlinear filter is derived. We also derive the nonlinear filtering equations associated with finitely-additive white noise formulation due to Kallianpur and Karandikar for the classical and quantum spin systems, and study existence and uniqueness of measure-valued solutions. 
\end{abstract}



\tableofcontents

\section{Introduction}
Classical and quantum spin systems are extensively studied in literature \cite{Ruelle1969, Simon1974, Simon1993, Bratteli1987, Bratteli1996}. In this paper we will take a stochastic partial differential equations approach as in \cite{DaPrato1995} for spin systems and develop a classical nonlinear filtering method for these models. For classical spin systems we will consider the diffusions on infinite dimensional torus studied by Holley and Stroock \cite{Holley1981} and also the Euclidean and lattice systems studied by Da Prato and Zabczyk \cite{DaPrato1995}. For quantum spin systems we consider the models studied by \cite{Albeverio2001} and \cite{DaPrato1995}. Nonlinear filtering for infinite dimensional problems related to nonlinear partial differential equations of fluid dynamics was initiated in \cite{Sritharan1995} and later developed for several fluid dynamic models in \cite{Sritharan2010, Sritharan2011, Fernando2013}. Nonlinear filtering for nonlinear stochastic partial differential equations for reaction diffusion type was originated in \cite{Hobbs1996}. In these papers the measure valued evolution equations of FKK type \cite{Fujisaki1972} and Zakai type\cite{Zakai1969} have been derived and existense and uniqueness of measure valued nonlinear filter was established. In this paper we will establish similar results for classical and quantum spin systems \cite{Albeverio2001,DaPrato1995} and then prove convergence and ergodicity type results similar to such results for finite dimensional nonlinear filtering problems initiated by Kunita\cite{Kunita1971} and further developed by several other authors \cite{Mazz1988, Stettner1989, Kunita1991, Ocone1996,Bhatt1995, Bhatt1999, Bhatt2000, Budhiraja2003}. We will also consider white noise nonlinear filtering equations developed by Kallianpur and Karandikar \cite{Kallianpur1985, Kallianpur1988, Kallianpur1991} and derive the corresponding equations for the spin systems and prove existence and uniqueness of measure-valued solutions. Other works on nonlinear filtering of classical spin systems with different perspective see \cite{Rebes2015}. We also note here that in this paper we develop classical nonlinear filtering for quantum spin systems even though white noise filtering in \textsection 6 has some attributes of quantum theory. Quantum nonlinear filtering for open quantum systems has been developed in a series of papers by V. P. Belavkin \cite{Belavkin1980, Belavkin1983, Belavkin1987, Belavkin1992a, Belavkin1992b} (see also the introductory exposition \cite{Bouten2007}) for the Hudson-Parthasarathy equation \cite{Hudson1984, Partha1992}. Belavkin type nonlinear filtering theory for unbounded Hamiltonian and noise operators in the Hudson-Parthasarathy equation is an open problem and will be addressed in the future.

\section{Unified Mathematical Formulation of Spin Systems -Stochastic Dissipative Systems and Invariant Measures}
In this section we will describe certain general results for infinite dimensional stochastic equations driven by cylindrical Wiener process as presented in \cite{DaPrato1995, DaPrato1996}. Let $H$ and $U$ be separable Hilbert spaces with norms $\|\cdot\|, \|\cdot\|_{U}$ and scalar products $\langle \cdot, \cdot \rangle$,  $\langle \cdot, \cdot \rangle_{U}$. The spaces of all bounded operators from $U$ in to $H$ and from $U$ in to $U$ will be denoted by ${\mathcal L}(U,H)$ and ${\mathcal L}(U)$ respectively. Let $(\Omega, {\mathcal F},{\mathcal F}_{t}, m)$ be a complete filtered probability space and $W(t), t\geq 0$ be an ${\mathcal F}_{t}$-adapted Wiener process defined on $\Omega$ with values in $U$ and with the covariance operator $Q\in {\mathcal L}(U)$. Hence for arbitary $\phi,\psi \in U$ and $t,s\geq 0$:
\begin{displaymath}
	E\left [ \langle W(t),\phi\rangle_{U}\langle W(s),\psi\rangle_{U}\right ]= t\wedge s \langle Q\phi, \psi\rangle_{U}. 
\end{displaymath}
We will consider the general stochastic equation
\begin{equation}
	dX=(AX+ F(X))dt +BdW, \label{eqn2.1}
\end{equation}
\begin{displaymath}
	X(0)=x.
\end{displaymath}
Here $B\in {\mathcal L}(U,H)$, $A\in {\mathcal L}(D(A);H)$ is a linear operator with $D(A)\subset H$, $F:D(F)\rightarrow H$ is a nonlinear operator with $D(F)\subset H$ and both $A$ and $F$ are dissipative in the sense defined below. An Important observation is that in the classes of problems studied in this paper the covariance operator is not of trace class: $\mbox{Tr}Q = \sum_{i}\langle Q \varphi_{i},\varphi_{i}\rangle_{U} =+\infty.$ In particular $Q=I$ is often explicitly invoked. 

Let $E, \|\cdot\|_{E}$ be a Banach space and $E^{*}$ its dual. For arbitrary $x\in E$ the subdifferential $\partial \|x\|_{E}$ of the norm $\|\cdot\|_{E}$ at $x$ is given by:
\begin{displaymath}
	\partial \|x\|_{E}:= \left \{ x^{*}\in E^{*}: \|x+y\|_{E}-\|x\|\geq x^{*}(y), \forall y\in E\right \}.	
\end{displaymath}
A mapping $G$ from $D(G)\subset E$ into $E$ is said to be dissipative in $E$ if for arbitrary $x,y\in D(G)$ there exists $z^{*} \in \partial \|x-y\|_{E}$ such that
\begin{equation}
	z^{*} (G(x)-G(y))\leq 0. \label{eqn2.2}
\end{equation}
If in addition, for some $\alpha>0$ the mapping $(I-\alpha G): D(G)\rightarrow E$ is surjective then $G$ is said to be $m$-dissipative.

Let $K\subset E$ is a Banach space embedded into $E$ then the part $G_{K}$. Then the part $G_{K}$ of $G$ in $K$ is defined as follows:
\begin{displaymath}
	G_{K}(x)=G(x) \mbox{ for } x\in D(G_{K}),
\end{displaymath}
where
\begin{displaymath}
	D(G_{K})=\left \{ x\in D(G)\cap K: G(x)\in K \right \}.
\end{displaymath}
We will now state the two essential hypothesis on operators $A, F, B$, the Wiener process $W(t), t\geq 0$ and on a Banach space $K, \|\cdot\|_{K}$ continuously and densely embedded into $H$.

\begin{Hyp} \label{H2.1}
	\begin{enumerate}
		\item There exists $\eta\in \mR$ such that the operators $A+\eta I$ and $F+\eta I$ are $m$-dissipative on $H$.
		\item The parts in $K$ of $A+\eta I$ and $F+\eta I$ are $m$-dissipative on $K$.
		\item $K\subset D(F)$ and $F$ maps bounded sets in $K$ into bounded sets in $H$.
		\item $K$ is a reflective space.
		\item $B\in {\mathcal L}(U,H)$.
	\end{enumerate}	
\end{Hyp}
Let $W_{A}(t), t\geq 0$ be the solution to the linear equation
\begin{equation}
	dZ(t)=AZ(t)dt + BdW(t), \label{eqn2.3}
\end{equation}
\begin{displaymath}
	Z(t)=W_{A}(t)=\int_{0}^{t} S(t-s)BdW(s), \hspace{.1in} t\geq 0,
\end{displaymath}
where $S(t), t\geq 0$, is the semigroup generated by $A$ on $H$.
\begin{Hyp}\label{H2.2}
	The process $W_{A}(t), t\geq 0$, is continuous in $H$, takes values in the domain $D(F_{K})$ of the part of $F$ in $K$ and for any $T>0$ we have 
	\begin{equation}
		\sup_{t\in [0,T]}\left (\|W_{A}(t)\|_{K} +\|F_{K}(W_{A}(t))\|_{K}\right )<+\infty, \hspace{.1in } m \mbox{ a.s. } \label{eqn2.4}
	\end{equation}
\end{Hyp}
\begin{Def}\label{Def2.1}
	An $H$-valued continuous ${\mathcal F}_{t}$-adapted process $X(t), t\geq 0$, is said to be a strong solution to \ref{eqn2.1} if it satisfies $m$-a.s. the equation
	\begin{equation}
		X(t)=x+ \int_{0}^{t} \left ( AX(s)+F(X(s))\right )	ds + BW(t), \hspace{.1in} t\geq 0, \label{eqn2.5}
	\end{equation}
	and it is a mild solution if it satisfies the following integral equation:
	\begin{equation}
		X(t)=S(t)x + \int_{0}^{t} S(t-s)F(X(s))ds +W_{A}(t), \hspace{.1in} t\geq 0. \label{eqn2.6}
	\end{equation}
	If for a $H$-valued process $X$, there exists a sequence $X_{n}$, of mild solutions to \ref{eqn2.1} such that $m$-a.s., $X_{n}\rightarrow X$ uniformly on any interval $[0,T]$, then $X$ is said to be a generalized solution to \ref{eqn2.1}. 
\end{Def}
Note that each strong solution is a mild solution and each mild solution is a generalized solution.

\begin{Thm}\label{Thm2.1}\cite{DaPrato1995, DaPrato1996}Assume that Hypothesis \ref{H2.1} and \ref{H2.2} are satisfied. Then for arbitrary $x\in K$ there exists a unique mild solution of \ref{eqn2.1} and for arbitrary $x\in H$ there exists a unique generalized solution $X(t,x), t\geq 0$ of \ref{eqn2.1}. If the operator $A$ and its part in $K$ are bounded then solutions for $x\in K$ are strong.
\end{Thm}
The above theorem implies that 
\begin{Thm}
$X(t,x),t\geq 0$ is a Markov process with transition group $P_{t}, t\geq 0$, given by
\begin{equation}
	P_{t}\varphi(x)=E[\varphi(X(t,x))], \hspace{.1in} t\geq 0, x\in H, \varphi \in C_{b}(H), \label{eqn2.7}
\end{equation}
where $C_{b}(H)$ denotes the space of all uniformly continuous and bounded functions on $H$. Moreover, the transition semigroup is Feller:
\begin{displaymath}
	P_{t}: C_{b}([0,\infty];H)\rightarrow C_{b}([0,\infty];H).
\end{displaymath}
\end{Thm}

The next theorem concerns with existence and uniqueness of invariant measures for \ref{eqn2.1} and asymptotic behavior of the transition semigroup $P_{t}, t\geq 0$.
\begin{Thm} \label{Thm2.2}\cite{DaPrato1995, DaPrato1996} If in addition to Hypotheses \ref{H2.1}, \ref{H2.2},
	\begin{enumerate}
		\item there exist $\eta_{1},\eta_{2}\in \mR$ such that $\omega=\eta_{1}+\eta_{2}>0$ and operators $A+\eta_{1}I$, $F+\eta_{2}I$ are dissipative in $H$,
		\item one has
		\begin{displaymath}
			\sup_{t\geq 0} E\left (\|W_{A}(t)\| +\|F(W_{A}(t))\|\right )<+\infty.
		\end{displaymath}
	\end{enumerate}
	Then there exists a unique invariant measure $\mu$ for the semigroup $P_{t}, t\geq 0$ such that $P^{*}_{t}\mu=\mu$. Moreover, for all bounded and Lipschitz continuous function $\varphi$ on $H$ one has
	\begin{equation}
		\vert P_{t}\varphi(x)-\int_{H}\varphi(y)\mu(y)\vert \leq (C+2 \|x\|) e^{-\omega t}\|\varphi\|_{\mbox{Lip}} \label{eqn2.8}
	\end{equation}
	where
	\begin{displaymath}
		C= \sup_{t\geq 0} E\left (\|W_{A}(t)\| +\frac{1}{\omega}\|F(W_{A}(t))\|\right ).
	\end{displaymath}
\end{Thm}
\section{Classical Spin Systems on Discrete Lattices}

\subsection{Diffusions on an Infinite Dimensional Torus}

We will discuss one of the simplest classical spin systems studied by Holley and Stroock \cite{Holley1981}. Let $(\Omega, \Sigma, m)$ be a complete probability space and let $x_{\gamma}(t)\in \mR^{\mZ^{d}}$ (resp.  $x_{\gamma}(t)\in \T^{\mZ^{d}}$,) $\gamma \in \mZ^{d}$ be a spin system that evolves according to:
\begin{equation}
x_{\gamma}(t,x)=x_{\gamma} +\int_{0}^{t}\sigma_{\gamma}(x(s,x))dB_{\gamma}(s)+\int_{0}^{t}b_{\gamma}(x(s,x))ds, \gamma \in \mZ^{d}, t>0, \label{eqn3.1}
\end{equation}
where for fixed $\gamma\in \mZ^{d}$, $B_{\gamma}$ is a one dimensional standard Brownian motion, $x \in \mR^{\mZ^{d}}	$ (resp. $x \in \T^{\mZ^{d}} )$, and $\sigma_{\gamma}(\cdot):\mR^{\mZ^{d}}\rightarrow \mR$, $b_{\gamma}(\cdot):\mR^{\mZ^{d}}\rightarrow \mR$ (resp. defined on $\T^{\mZ^{d}}$) and their derivatives are smooth bounded functions such that $\sigma_{\gamma}(x)=\sigma_{\gamma}(y)$ if $x_{l}=y_{l}$ for all $l\in Z^{d}$ satisfying $\vert l-\gamma\vert \leq L$. We will start with a path-wise solvability:

\begin{Lem}\cite{Holley1981} \label{Lem3.1} For each $x\in \mR^{\mZ^{d}}$ there exists a unique solution $x_{\gamma}(\cdot,x)$ to \ref{eqn3.1}. Moreover, if $x,y\in \mR^{\mZ^{d}}$, then for each $T>0$:
	\begin{displaymath}
		E\left [\sum_{\gamma\in \mZ^{d}}\sup_{0\leq t\leq T}\frac{1}{2^{\vert \gamma\vert}}\vert x_{\gamma}(t,x)-x_{\gamma}(t,y)\vert^{2}\right]\leq A_{T}\left (\sum_{\gamma\in \mZ^{d}}\frac{1}{2^{\vert \gamma\vert}}\vert x_{\gamma}-y_{\gamma}\vert^{2}\right ),
	\end{displaymath}
where $A_{T}<\infty$ depends only on $T>0, L,$ and $\max_{\gamma}\|\sigma_{\gamma}\|_{C^{1}(\mR^{\mZ^{d}})}$. Finally, if $x^{N}(\cdot,x)_{\gamma}$ is defined by \ref{eqn3.1} with $\sigma_{\gamma}(\cdot)=b_{\gamma}(\cdot)\equiv 0$ for $\vert \gamma \vert >N$ then for all $T>0$, 
\begin{displaymath}
E\left [ \sum_{\gamma \in \mZ^{d}} \sup_{0\leq t\leq T}\frac{1}{2^{\vert \gamma\vert }}\vert x_{\gamma}(t,x)-x^{N}_{\gamma}(t,x)\vert^{2}\right ] \rightarrow 0 \mbox{  as  } N\rightarrow \infty.
\end{displaymath}
\end{Lem}
Let us now recall the solvability of martingale problem given by Holley and Stroock \cite{Holley1981}.
Let $\T=\left \{z\in \mC; \vert z\vert =1\right \}$. For points in $\T^{\mZ}$ we will use $\eta$ to denote both the point in $\T^{\mZ}$ as well as the element $\alpha$ of $([0,2\pi))^{\mZ}$ such that $\eta_{k}=e^{i\alpha_{k}}, k\in \mZ$. Let ${\mathcal D}$ denote the set of smooth functions on $\T^{\mZ}$ (sometimes called cylindrical test  functions) which depend only on a finite number of coordinates.
\begin{Thm} \cite{Holley1981}
Let $\Omega$ denote the Polish space $C([0,\infty);\T^{\mZ})$. For $\omega \in \Omega, \eta(t,\omega)\in \T^{\mZ}$ is the position of $\omega$ at time $t\geq 0$.	We set $\mcF_{t}:=\sigma\left(\eta(s); 0\leq s\leq t \right)$ to be the cannonical filtration and $\mcF=\sigma \left(\cup_{t\geq 0}\mcF_{t}\right)$ and then $\mcF$ coincides with the Borel field over $\Omega$. Let $\sigma_{k},b_{k}, k\in \mZ$ be smooth functions of $\T^{\mZ}$ such that for a given $k\in \mZ$ the functions $\sigma_{k}$ and $b_{k}$ depend only on $\left \{ \eta_{l};\vert l-k\vert \leq L\right \}$. Also assume that $\sigma_{k}$ and $b_{k}$  and each of their derivatives are bounded independent of $k\in \mZ$. For $f\in {\mathcal D}$, define
\begin{displaymath}
	{\mathcal L}f(\eta) = \sum_{k\in \mZ}\left ( \frac{1}{2} \sigma_{k}^{2}(\eta)\frac{\partial^{2} f}{\partial \eta_{k}^{2}}(\eta) +b_{k}(\eta) \frac{\partial f}{\partial \eta_{k}}(\eta)  \right ), \hspace{.1in} \eta \in \T^{\mZ}.
\end{displaymath}
Then for each $\eta\in \T^{\mZ}$ there is exactly one probability measure $P_{\eta}$ on $(\Omega, \mcF)$ such that $P_{\eta}(\eta(0)=\eta)=1$ and $\left \{f(\eta(t))-\int_{0}^{t}{\mathcal L}f(\eta(s))ds, {\mathcal F}_{t}, P_{\eta}\right \}$ is a martingale for each $f\in {\mathcal D}$. Moreover, if $\Phi:\mR\rightarrow \mZ$ is the map defined so that $\Phi_{k}(x)=\eta_{k}$, where $\eta_{k}\in [0,2\pi) $ and $x_{k}=\eta_{k}$ mod $(2\pi)$, then for any $n\geq 1$, $F\in B((\T^\mZ)^{n})$, and $0\leq t_{1}<t_{1} <\cdots <t_{n}$:
\begin{displaymath}
	E^{P_{\eta}}\left [ F(\eta(t_{1}),\cdots, \eta(t_{n}))  \right]=E^{m}\left[F\circ \Phi^{n}(x(t_{1},\eta),\cdots, x(t_{n},\eta)\right],
\end{displaymath}
where $\Phi^{n}=\Phi \otimes \cdots \otimes \Phi$ ($n$-times ), and  $x(\cdot,\eta)$ 
 is the solution to \ref{eqn3.1} with $x=\eta$ and the $\sigma_{k}$ and $b_{k}$ in \ref{eqn3.1} replaced by $\sigma_{k}\circ \Phi$ and $b_{k}\circ \Phi$, respectively. Finally, the family $\left \{P_{\eta}; \eta \in \T^{\mZ}\right \}$ is Feller continuous and strong Markov.
\end{Thm}	

\subsection{Unbounded Classical Spin Systems}

Unbounded classical spin systems have been studied by several authors to establish construction of Markov semigroup, existence and uniqueness of invariant measures, ergodicity and exponential convergence of Markov semigroups to equilibrium state characterized by the invariant measure \cite{Zagarlinski1996,DaPrato1995,DaPrato1996}. We will recall some of the essential results below. We consider the Markov process $X=\{X_{\gamma}\}_{\gamma\in \mZ^{d}}$ satisfying an infinite system of Ito's equations
\begin{equation}
dX_{\gamma}(t) =\left (\sum_{j}a_{\gamma j}X_{j}(t) +f(X_{\gamma}(t))\right)dt +dW_{\gamma}(t), \label{eqn3.2}
\end{equation}
\begin{displaymath}
	X_{\gamma}(0)=x_{\gamma}, \hspace{.1in} \gamma\in \mZ^{d}, t\geq 0,
\end{displaymath}
where $W_{\gamma}$, $\gamma\in \mZ^{d}$ is a Wiener process on $(\Omega, \Sigma, m)$ with values in $U=l^{2}(\mR^{d})$ and with the covariance operator $Q\in {\mathcal L}(U)$. In particular, if $Q=I$ then the processes $W_{\gamma},\gamma\in \mZ^{d}$, are independent standard real valued Wiener processes. An interesting special case of the above system is:
\begin{equation}
	dX_{\gamma}(t) =\left ( (\Delta_{d}-\alpha)X_{\gamma}(t) +f(X_{\gamma}(t))\right)dt +dW_{\gamma}(t), \label{eqn3.3}
\end{equation}
\begin{displaymath}
	X_{\gamma}(0)=x_{\gamma}, \hspace{.1in} \gamma\in \mZ^{d}, t\geq 0,
\end{displaymath}
where $\Delta_{d}$ is the $d$-dimensional discrete Laplacian and $\alpha$ is a constant.

In the light of general theory presented in Section 2, we define operators $A$ and $F$ as follows:
\begin{displaymath}
	A(X_{\gamma})=\left (\sum_{j\in \mZ^{d}} a_{j\gamma}X_{j}\right ), \hspace{.1in} X=(X_{\gamma})\in H,
\end{displaymath}
\begin{displaymath}
	F(X_{\gamma})=(f(X_{\gamma})), \hspace{.1in} X=(X_{\gamma})\in H.
\end{displaymath}
The following result \cite{DaPrato1995, DaPrato1996, Peszat2007} is useful in establishing the boundedness of $A$.
\begin{Thm} \label{Thm3.2} Assume that
	\begin{enumerate}
		\item $\sup_{\gamma \in \mZ^{d}}\sum_{j\in \mZ^{d}} \vert a_{\gamma j}\vert =\alpha <+\infty.$
		\item There exists $\beta >0$ such that
		\begin{displaymath}
			\sum_{\gamma \in \mZ^{d}}\vert a_{\gamma j} \vert \rho(\gamma)\leq \beta \rho (j),  \hspace{.1in} j\in \mZ^{d}.
		\end{displaymath}
	Then the formula 
	\begin{displaymath}
	A_{p}(X_{\gamma})=\left (\sum_{j\in \mZ^{d}} a_{j\gamma}X_{j}\right ), \hspace{.1in} X=(X_{\gamma})\in l^{p}_{\rho}(\mZ^{d}),	
	\end{displaymath}
defines a bounded operator on $l^{p}_{\rho}(\mZ^{d})$, (i.e. $A\in {\mathcal L}(l^{p}_{\rho},l^{p}_{\rho})$) for all $p\in [1,+\infty]$, with the norm not greater than
\begin{displaymath}
	\alpha^{1/q}\beta^{1/p}, \hspace{.1in} \frac{1}{p}+\frac{1}{q}=1.
\end{displaymath}
	\end{enumerate}
In particular $A\in {\mathcal L}(l^{2}_{\rho},l^{2}_{\rho})$, with operator norm less than or equal to $\sqrt{\alpha \beta }$.
\end{Thm}
\begin{Cor} \label{Cor3.1}
	If, for some $\alpha<\infty$,
	\begin{displaymath} 
		\sup_{k\in \mZ^{d}}\sum_{j\in \mZ^{d}}\vert a_{k,j}\vert \leq \alpha \mbox{ and } \sup_{j\in \mZ^{d}}\sum_{k\in \mZ^{d}}\vert a_{k,j}\vert \leq \alpha
	\end{displaymath}
	then 	$A\in {\mathcal L}(l^{p},l^{p})$ with norm $\leq \alpha$.
\end{Cor}	
For example, the matrix coefficients $a_{\gamma, j}$ and positive weight $\rho$ can be assumed to satisfy the following conditions:
\begin{Def} \label{Def3.1}
	The interactions $a_{\gamma,j}$ have {\it bounded range} if there are constants $R,M>0$ such that
	\begin{equation}
		a_{\gamma,j}=0 \mbox{ if } \vert \gamma -j \vert > R \mbox{ and } 
			\vert a_{\gamma,j}\vert \leq M \mbox{ for all  } \gamma, j \in \mZ^{d}. \label{eqn3.4}
	\end{equation}
\end{Def}

\begin{Lem} \label{Lem3.2}
	Assume {\it bounded range } and that 
	\begin{equation} 
		\vert \frac{\rho(\gamma)}{\rho(j)}\vert \leq M \mbox{ if } \vert \gamma -j \vert \leq R,\mbox{ and }\sum _{\gamma \in \mZ^{d}}\rho(\gamma) <+\infty, \label{eqn3.5}
	\end{equation}
then $A\in {\mathcal L}(l^{p}_{\rho},l^{p}_{\rho})$ for every $p\geq 1$.
\end{Lem}

Examples of weights $\rho$ that fit the conditions of 
 Lemma \ref{Lem3.2} are:
\begin{displaymath} 
		\rho_{\kappa}(k)=e^{-\kappa \vert k\vert } \mbox{ or } \rho^{\kappa}(k)=\frac{1}{1+\kappa \vert k\vert^{r}}, \kappa >0, r>d, k\in \mZ^{d}.
		\end{displaymath}
	\begin{Def} \label{Def3.2}
	The local interaction function $f:\mR\rightarrow \mR$ is such that $f=f_{0}+f_{1}$ where     $f_{1}$ is Lipchitz continuous and $f_{0}+\eta \zeta, \zeta \in \mR$ is continuous and decreasing for some $\eta \in \mR$ and for some  $s\geq 1$ and $c_{0}>0$,
	\begin{displaymath}
		\vert f_{0}(\zeta)\vert \leq c_{0}( 1 + \vert \zeta \vert^{s}), \zeta \in \mR.
	\end{displaymath}
\end{Def}
	For example with  $s=2n+1$ and 
	\begin{displaymath}
		f_{0}(\zeta) = -\zeta^{2n+1}+\sum_{k=0}^{2n}b_{k}\zeta^{2n-k}, \hspace{.1in} \zeta \in \mR, n \in \mN\cup \{0\}.
	\end{displaymath}
	
\begin{Thm}\label{Thm3.3}
	Let $ H=l_{\rho}^{2}(\mZ^{d})$, $K=l_{\rho}^{2s}(\mZ^{d})$ and let the operators $A$ and $F$ be as above. Then
	\begin{enumerate} 
			\item For arbitrary $x\in K=l_{\rho}^{2s}(\mZ^{d})$ there exists a unique strong solution $X(t,x), t\geq 0$. 
			\item For arbitrary $x \in H=l_{\rho}^{2}(\mZ^{d})$, there exists a unique generalized solution $X(t,x), t\geq 0$, and the transition semigroup 
			\begin{displaymath} 
				P_{t}\varphi(x)= E(\varphi (X(t,x))), \hspace{.1in} t\geq 0, x \in H, \varphi \in C_{b}(H),
			\end{displaymath}
			is Feller.
	\end{enumerate}
\end{Thm}

\begin{Thm}\label{Thm3.4}
	Assume that conditions \ref{eqn3.4}, \ref{eqn3.5} are satisfied. Let $ H=l_{\rho}^{2}(\mZ^{d})$, $K=l_{\rho}^{2s}(\mZ^{d})$ and let the operators $A$ and $F$ be given as above. 	Then\begin{enumerate}
		\item For arbitrary $x\in K=l_{\rho}^{2s}(\mZ^{d})$ there exists a unique strong solution $X(t,x), t\geq 0$ of \ref{eqn3.2}.
		\item For arbitrary $x \in H=l_{\rho}^{2}(\mZ^{d})$, there exists a unique generalized solution $X(t,x), t\geq 0$ of \ref{eqn3.2}  and the transition semigroup 
		\begin{displaymath}
			P_{t}\varphi(x)= E(\varphi (X(t,x))), \hspace{.1in} t\geq 0, x \in H, \varphi \in C_{b}(H),
		\end{displaymath}
	is Feller.
	\end{enumerate}
\end{Thm}

\begin{Thm}\label{Thm3.5}
	Assume, in addition to conditions of \ref{Thm3.4}
	that for an $\eta >0$, operator $A+I\eta$, restricted to $l^{2}(\mZ^{d})$ is dissipative, that $f_{0}$ is decreasing and $\eta-\|f_{1}\|_{\mbox{Lip}}>\omega > 0$. Then there exists $\kappa_{0}>0$ such that in the spaces $l_{\rho}^{2}(\mZ^{d})$ with $\rho=\rho^{\kappa}$  or with $\rho=\rho_{\kappa}$, given as above, and $\kappa \in ]0,\kappa_{0}[$, equation \ref{eqn3.2} has unique generalized solutions. The semigroup $P_{t}, t\geq 0$ has a unique invariant measure $\mu$ on $H$ and there exists $C>0$ such that, for any bounded and Lipschitz function $\varphi$ on $H$, all $t>0$ and all $x\in H$:
	\begin{displaymath}
	\vert P_{t}\varphi(x)-\int_{H}\varphi(x)\mu(dx)\vert \leq (C+2\vert x\vert ) e^{-\omega t}\|\varphi\|_{\mbox{ Lip}}.
	\end{displaymath}
\end{Thm}
These results have been extended to spin systems forced by L\'evy noise in \cite{Peszat2007}.

\subsection{Euclidean Field Theory and Spin Systems on Euclidean Spaces }
Spin systems on Euclidean spaces can be modeled by the following stochastic partial differential equation:
\begin{equation}
dX(t,\zeta)=[(\Delta-\alpha)X(t,\zeta) + f(X(t,\zeta))]dt+dW(t,\zeta), \label{eqn3.6}
\end{equation}
\begin{displaymath}
X(0,\zeta)=x(\zeta), \hspace{.1in} \zeta \in \mR^{d}, t>0,
\end{displaymath}
where $\Delta$ is the Laplace operator, $\alpha$ is a nonnegative constant and $W(t,\zeta), \zeta \in \mR^{d}$, $t>0$, an infinite dimensional Wiener process with covariance operator $Q$, defined on a probability space $(\Omega, {\mathcal F}, m)$.
\begin{Hyp} \label{H3.1}
	\begin{enumerate} 
					\item 
						\begin{equation}
	Q=I \mbox{  on } U=L^{2}(\mR^{d}), \label{eqn3.7}
			\end{equation}
	or
			\item $Q$ is a convolution operator on $U=L^{2}(\mR^{d})$:\\ 
			\begin{equation}	
			Qu(\zeta)=q\ast u(\zeta), \zeta\in \mR^{d}, u\in L^{2}(\mR^{d}), \label{eqn3.8}
		\end{equation}
	where $ q(\zeta)=\int_{\mR^{d}}e^{i \zeta\cdot \eta}g(\eta)d\eta, \zeta\in \mR^{d},$ 	with $g\geq 0$ and $q, g\in L^{1}(\mR^{d})$.
	\end{enumerate}
\end{Hyp}
\begin{Thm}\label{T3.4}
	Assume that $f$ is as in Definition \ref{Def3.2} and either
	\begin{enumerate}
		\item $d=1$ and \ref{eqn3.7} or \ref{eqn3.8} holds or
		\item $d>1$ and \ref{eqn3.8} holds.
	\end{enumerate}
Then the equation \ref{eqn3.6} has a unique generalized solution in $L^{2}_{\rho}(\mR^{d})$ where $\rho$ is given either by $\rho_{\kappa}$  or by $\rho^{\kappa}$. If $x\in L^{2s}_{\rho}(\mR^{d})$ then the generalized solution is mild.
\end{Thm}

\begin{Thm}
In addition to conditions of Theorem \ref{T3.4} assume that the function $f_{0}$ is decreasing and $\alpha -\|f_{1}\|_{\mbox{ Lip }} >\omega >0$. Then there exists $\kappa_{0}>0$ such that the semigroup $P_{t}, t\geq 0$ corresponding to solution of equation \ref{eqn3.6} has a unique invariant measure both in $H=L^{2}_{\rho^{\kappa}	}(\mR^{d})$ and $H=L^{2}_{\rho_{\kappa}	}(\mR^{d})$, for any $\kappa \in ]0,\kappa_{0}[.$ Moreover there exists $C>0$ such that for any bounded Lipschitz function $\varphi$ on $H$, all $t>0$ and all $x\in H$
\begin{equation}
	\vert P_{t}\varphi(x)-\int_{H}\varphi(x)\mu(dx)\vert \leq (C+2) \|x\|e^{-\omega t}\|\varphi\|_{\mbox{Lip }}.
\end{equation}
\end{Thm}

\section{Quantum Lattice Systems}
One of the simplest quantum lattice system would be a combination of the classical and Euclidean spin systems formulated in the previous sections (see \cite{DaPrato1995,DaPrato1996}):
\begin{equation}
dX_{\gamma}(t,\zeta)=\left ( {\mcA}X_{\gamma}(t,\zeta) +\sum_{j\in\mZ^{d}}a_{\gamma j} X_{\gamma}(t,\zeta)+ {\mcF}(X_{\gamma}(t,\zeta))\right )dt +dW_{\gamma}(t,\zeta),	\label{eqn4.1}
\end{equation}
\begin{displaymath}
	X_{\gamma}(0,\zeta)=x_{\gamma}(\zeta), \hspace{.1in} \zeta \in [0,1], t>0,
\end{displaymath}
where $\mcA$, $\mcF$ are in general unbounded, respectively linear and nonlinear operators on a Hilbert space ${\mathcal H}$, $(a_{\gamma j})_{\gamma,j\in \mZ^{d}}$ is a given matrix with real elements and $(W_{\gamma})_{\gamma \in \mZ^{d}}$ is a family of independent cylindrical Wiener processes on ${\mathcal H}$.
We will assume for simplicity that ${\mathcal H}=L^{2}(0,1)$, ${\mcA}=\frac{d^{2}}{d\zeta^{2}}-\alpha$,
\begin{displaymath}
D({\mcA})=\left \{ x\in H^{2}(0,1): x(0)=x(1), x'(0)=x'(1)\right \},
\end{displaymath}
\begin{displaymath}
	{\mcF}(x)(\zeta)=f(x(\zeta)), \hspace{.1in} \zeta \in [0,1],
\end{displaymath}
\begin{displaymath}
	D({\mcF})=\left \{x\in L^{2}(0,1): f(x)\in L^{2}(0,1)\right \},
\end{displaymath}
and $f$ and the matrix $\{a_{\gamma j}\}_{\gamma, j \in \mZ^{d}}$ satisfy conditions of Definition \ref{Def3.2} and Definition \ref{Def3.1} respectively. In order to frame \ref{eqn4.1} in to the general form of equation \ref{eqn2.1} we define:
\begin{displaymath}
	H=l^{2}_{\rho}(L^{2}(0,1))=\left \{(x_{\gamma})\in {\mathcal H}^{(\mZ^{d})}: \sum_{\gamma \in \mZ^{d}} \rho(\gamma)\|x_{\gamma}\|^{2}_{\mathcal H}<\infty \right \},
\end{displaymath}
\begin{displaymath}
	K=l^{2s}_{\rho}(L^{2s}(0,1))=\left \{(x_{\gamma})\in (L^{2s}(0,1))^{(\mZ^{d})}: \sum_{\gamma \in \mZ^{d}} \rho(\gamma)\|x_{\gamma}\|^{2}_{L^{2s}(0,1)}<\infty \right \},
\end{displaymath}
where $\rho=\rho^{\kappa}$ or $\rho=\rho_{\kappa}$, $\kappa>0$ as in Section 3.2. Let $A=A_{0}+ A_{1}$ where $A_{1}$ is a bounded linear operator on $H$ given by
\begin{displaymath}
	A_{1}(x_{\gamma})= \left (\sum_{j\in\mZ^{d}}a_{\gamma j} X_{\gamma}\right), x\in D(A_{1})=H,
\end{displaymath}
and $A_{0}(x_{\gamma})=({\mcA}x_{\gamma}),$ $x=(x_{\gamma})\in D(A_{0})$,
\begin{displaymath}
	D(A_{0})=\left \{ (x_{\gamma})\in H: \sum_{\gamma \mZ^{d}} \rho(\gamma)\|{\mcA}x_{\gamma}\|^{2}_{\mathcal H}<\infty \right\}.
\end{displaymath}
$A_{1}$ is bounded by a generalization of Theorem \ref{Thm3.2}. and $A_{0}+\eta$ on $H$ and its restriction $A_{0\rho}$ to $K$ are $m$-dissipative for sufficiently small $\eta$. Let us define
\begin{displaymath}
F(x_{\gamma})=({\mcF}x_{\gamma}), \hspace{.1in} x=(x_{\gamma})\in D(F),
\end{displaymath}
\begin{displaymath}
	D(F) =\left \{ (x_{\gamma})\in H: \sum_{\gamma \in \mZ^{d}} \rho(\gamma)\|{\mcF}x_{\gamma}\|^{2}_{\mathcal H}<\infty \right\}.
\end{displaymath}

\begin{Thm}\label{T4.1}
Assume that conditions of Definition \ref{Def3.2} hold. Then for arbitrary $x\in H$ equation \ref{eqn4.1} has a unique generalized solution $X(\cdot,x)$. If $x\in K$ the solution is mild.
\end{Thm}
\begin{Thm}
	In addition to conditions of Theorem \ref{T4.1} assume that the function $f_{0}$ is decreasing and $\alpha -\|f_{1}\|_{\mbox{ Lip }} >\omega >0$. Then there exists $\kappa_{0}>0$ such that the semigroup $P_{t}, t\geq 0$ corresponding to solution of equation \ref{eqn4.1} has a unique invariant measure $\mu$ both in $H=L^{2}_{\rho^{\kappa}	}(\mR^{d})$ and $H=L^{2}_{\rho_{\kappa}	}(\mR^{d})$, for any $\kappa \in ]0,\kappa_{0}[$. Moreover there exists $C>0$ such that for any bounded Lipschitz function $\varphi$ on $H$, all $t>0$ and all $x\in H$
	\begin{equation}
		\vert P_{t}\varphi(x)-\int_{H}\varphi(x)\mu(dx)\vert \leq (C+2) \|x\|e^{-\omega t}\|\varphi\|_{\mbox{Lip }}.
	\end{equation}
\end{Thm}

\section{Nonlinear Filtering Formulation and Filtering Equations}
\subsection{Stochastic Calculus Method of Nonlinear Filtering}

Nonlinear filtering theory for nonlinear stochastic partial differential equations of fluid mechanics and reacting and diffusing systems was initiated in \cite{Sritharan1995, Hobbs1996} and further developed for stochastic models of fluid mechanics with L\'evy noise in \cite{Sritharan2010, Sritharan2011, Fernando2013}. The nonlinear filtering problem for the class of spin systems studied in this paper is formulated as follows. The key mathematical equations of  stochastic calculus method of nonlinear filtering are the Kallianpur-Striebel formula, the Fujisaki-Kallianpur-Kunita (FKK) equation, the Zakai equation and the Kunita's semigroup versions of FKK and the Zakai equation, all will be developed below. We will also establish the equivalency of FKK, Zakai with their Kunita counter parts in the spirit of Szpirglas \cite{Szpirglas1978}. We call the process $X(t)\in H,t\geq 0$ the signal process and it is governed by the evolution formulated in the earlier sections as:
\begin{equation}
dX=(AX+ F(X))dt +BdW, \label{eqn5.1}
\end{equation}
\begin{displaymath}
	X(0)=x.
\end{displaymath}
We note here that $H$ is compact for the finitely truncated classical spin system, locally compact for the classical spin system and general Hilbert space for the Euclidean fields and also quantum spin lattice system as formulated in Section 2, 3 and 4. \\

The measurement process $Y(t)\in \mR^{N}, t\geq 0$ is defined as
\begin{equation}
	dY=h(X)dt + dZ, \label{eqn5.2}
\end{equation}
where $Z(t)\in \mR^{N}$ is a $N$-dimensional Wiener process which may or may not be correlated with $W(t)$. Let us denote by ${\mathcal F}^{Y}_{t}$ the filtration generated by the sensor data over the time period $0\leq s\leq t$ (sigma algebra generated by the back measurement):
\begin{displaymath}
	{\mathcal F}^{Y}_{t}= \sigma \left \{Y_{s}, 0\leq s\leq t\right \}.
\end{displaymath}
Our goal is to study the time evolution of the conditional expectation $E[f(X(t))\vert {\mathcal F}^{Y}_{t}], t\geq 0$ called the nonlinear filter, where $f$ is some measurable function. It is also the least square best estimate of $f(X(t))$ given the back measurements. We will first consider a special case where the signal (spin system dynamics) noise and observation noise are independent. The following Bayes formula known as the Kallianpur-Striebel formula can be derived following \cite{Kallianpur1969}: Let $g$ be ${\mathcal F}^{X}_{t}$-measurable function $g(\cdot):H\rightarrow \mR$ integrable on $(\Omega, \Sigma, m)$, $0\leq t\leq T$. We may assume that processes in \ref{eqn5.2} be defined on the product space $(\Omega\times \Omega_{0}, \Sigma\times \Sigma_{0},m\times m_{0})$ where $Z$ is a Wiener process on $(\Omega_{0},\Sigma_{0}, m_{0})$. Then
\begin{equation}
	E_{m}[g\vert {\mathcal F}^{Y}_{t}] =\frac{\vartheta^{Y}(g)}{\vartheta^{Y}(1)} \label{eqn5.3}
\end{equation}
where 
\begin{displaymath}
	\vartheta^{Y}_{t}(g)=\int_{\Omega} g q_{t}dm
\end{displaymath}
with
\begin{displaymath}
	q_{t}=\exp \left \{ \int_{0}^{t} h(X(s))\cdot d Y(s) -\frac{1}{2} \int_{0}^{t}\vert h(X(s))\vert^{2}ds\right \}.
\end{displaymath}
From the Kallianpur-Striebel formula we can derive the nonlinear filtering equations of Zakai and FKK type using It\'o formula for the special case of signal and noise independent. However we can proceed as follows to get the filtering equations for the general case.

Let us define the formal infinitesimal generator ${\mathcal A}$ of the process $X(t)$ using the transition semigroup $P_{t}$ of $X(t)$ in the following way: for a function $f(\cdot):H\rightarrow \mR$, we set $P_{t}f(x)=E[f(X(t))\vert X(0)=x]$. 

Let us consider the measure space $(H, {\mathcal B}(H)) $ where ${\mathcal B}(H)$ is the Borel algebra generated by closed/open subsets of $H$. Let ${\mathcal B}_{b}(H)$ be the class of bounded ${\mathcal B}(H)$-measurable functions on $H$. For $f_{n},f \in {\mathcal B}_{b}(H)$ we say that $f_{n}\rightarrow f$ weakly if $f_{n}\rightarrow f$ pointwise and the sequence $f_{n}$ is uniformly bounded.

Let ${\mathcal D}({\mathcal A})$ be the class of functions $f$ in ${\mathcal B}_{b}(H)$ such that there exists $f_{0}\in {\mathcal B}_{b}(H)$ satisfying
\begin{equation}
	(P_{t}f)(x)=f(x) + \int_{0}^{t}(P_{s}f_{0})(x)ds, \hspace{.1in} \forall x\in H.
\end{equation}
Note that $f_{0}$ is uniquely determined by the above equation. Hence we set ${\mathcal A}f=f_{0}$ and formally define the extended generator ${\mathcal A}$ as:
\begin{displaymath}
	{\mathcal A}f(x):= \mbox{  weak-limit}_{t\downarrow 0^{+}}\frac{P_{t}f(x)-f(x)}{t}, \hspace{.1in} f\in D({\mathcal A}),
\end{displaymath}
where  
\begin{displaymath}
 D({\mathcal A})= \left \{f\in {\mathcal B}_{b}(H): \mbox{  weak-limit }_{t\downarrow 0^{+}}\frac{P_{t}f(x)-f(x)}{t} = {\mathcal A}f(x) \mbox{  exists} \right\}.	
\end{displaymath}
Since the domain $D({\mathcal A})$ cannot easily be defined explicitly, we will also work with a sub-class of functions known as cylindrical test functions ${\mathcal D}_{\mbox{cyl}}$ defined as follows:
\begin{displaymath}
{\mathcal D}_{\mbox{cyl}} :=\left \{f \in C_{b}(H): f(X)=\varphi (\langle X, e_{1}\rangle_{H},\cdots, \langle X, e_{N}\rangle_{H} ), e_{i}\in D(A), i=1,\cdots, N, \varphi\in C^{\infty}_{0}(\mR^{N})\right \}.
\end{displaymath}
Then $D_{x}f(x)\in D(A)$ and also ${\mathcal D}_{\mbox{cyl}}\subset D({\mathcal A})$. The formal generator of the stochastic process $X(t)$ is given by
\begin{equation}
	{\mathcal A}\Phi(x)=\frac{1}{2}\mbox{Tr}(B^{*}QBD^{2}_{x}\Phi(x))+\langle Ax+F(x),D_{x}\Phi \rangle, \hspace{.1in} \Phi\in {\mathcal D}_{\mbox{cyl}}. \label{eqn5.4}
\end{equation}

We now note that the solvability theorems established in the previous section imply by It\'o's formula the following result.
\begin{Pro}
Define 
\begin{equation}
M_{t}(f)=f(X(t))-f(x)-\int_{0}^{t}{\mathcal A}f(X(s))ds, \hspace{.1in} f\in {\mathcal D}_{\mbox{cyl}}, \label{eqn5.5}
\end{equation}
where $X(t)$ is the spin process defined by the stochastic partial differential equation \ref{eqn5.1}. Then $(M_{t}(f), {\mathcal F}_{t}, m)$ is a martingale.
\end{Pro}
In fact applying It\'o formula \cite{DaPrato2014} to the stochastic process defined by \ref{eqn5.1}
\begin{displaymath}
f(X(t))=f(x)+\int_{0}^{t}{\mathcal A}f(X(s))ds	+\int_{0}^{t}(\frac{\partial f}{\partial x}(X(s), BdW(s)), \hspace{.1in} f\in {\mathcal D}_{\mbox{cyl}}.
\end{displaymath}
Let us now state a general martingale lemma and then specialize it to the filtering problem.
\begin{Lem} Let ${\mathcal F}_{t}$ and ${\mathcal G}_{t}$ be filtrations with ${\mathcal G}_{t}\subset {\mathcal F}_{t}$. Suppose that 
	\begin{displaymath}
		M_{t}^{\mathcal F}=U(t)-U(0)-\int_{0}^{t}V(s)ds
	\end{displaymath}
is an ${\mathcal F}_{t}$-martingale. Then
\begin{displaymath}
M^{\mathcal G}_{t}=	E[U(t)\vert {\mathcal G}_{t}]-E[U(0)\vert {\mathcal G}_{0}]-\int_{0}^{t}E[V(s)\vert {\mathcal G}_{s}]ds
\end{displaymath}
is a ${\mathcal G}_{t}$-martingale.
\end{Lem}
We now state this property in the nonlinear filtering context:
\begin{Thm}
Define
\begin{equation}
	M_{t}^{Y}(f)=E[f(X(t))\vert {\mathcal F}_{t}^{Y}]-f(x) -\int_{0}^{t}E[{\mathcal A}f(X(s))\vert {\mathcal F}_{s}^{Y}]ds, \hspace{.1in} f\in {\mathcal D}_{\mbox{cyl}}.  \label{eqn5.6}
\end{equation}	
Then $(M_{t}^{Y}(f), {\mathcal F}_{t}^{Y}, m)$ is a martingale.
\end{Thm}

{\bf Proof:}
First note that
\begin{displaymath}
E\left [M^{Y}_{f}(f)\vert {\mathcal F}^{Y}_{s}\right ]-M^{Y}_{s}(f)=E\left [M^{Y}_{f}(f)-M^{Y}_{s}(f)\vert {\mathcal F}^{Y}_{s}\right ].	
\end{displaymath}
Noting that ${\mathcal F}_{t}^{Y}\subset {\mathcal F}_{t}$, we consider, for $s<t$,
\begin{displaymath}
	E\left [M^{Y}_{f}(f)-M^{Y}_{s}(f)\vert {\mathcal F}^{Y}_{s}\right ]
\end{displaymath}
\begin{displaymath}
	= E\left [E[f(X(t))\vert {\mathcal F}_{t}^{Y}]-E[f(X(s))\vert {\mathcal F}_{s}^{Y}]\vert {\mathcal F}^{Y}_{s}\right ]- E\left [\int_{s}^{t}E[{\mathcal A}f(X(r))\vert {\mathcal F}_{r}^{Y}]dr\vert {\mathcal F}^{Y}_{s}\right ].
\end{displaymath}
Using the tower property of the conditional expectation we conclude that
\begin{displaymath}
	E\left [M^{Y}_{f}(f)-M^{Y}_{s}(f)\vert {\mathcal F}^{Y}_{s}\right ]= E[f(X(t))\vert {\mathcal F}_{s}^{Y}]-E[f(X(s))\vert {\mathcal F}_{s}^{Y}]-  \int_{s}^{t}E[{\mathcal A}f(X(r))\vert {\mathcal F}_{s}^{Y}]dr
\end{displaymath}
\begin{displaymath}
	=E\left [ f(X(t))-f(X(s))-\int_{s}^{t}{\mathcal A}f(X(r))dr \vert {\mathcal F}_{s}^{Y}\right ]
\end{displaymath}
\begin{displaymath}
= E\left [M_{t}(f)-M_{s}(f)\vert {\mathcal F}^{Y}_{s}\right]= E\left [E\left [M_{t}(f)-M_{s}\vert{\mathcal F}_{s}\right] \vert {\mathcal F}^{Y}_{s}\right ]=0,
\end{displaymath}
since $M_{t}(f)$ is a $({\mathcal F}_{t},m)$- martingale. Hence $E\left [M^{Y}_{f}(f)\vert {\mathcal F}^{Y}_{s}\right ]= M^{Y}_{s}(f)$.

Using Fujisaki-Kallianpur-Kunita (1972) \cite{Fujisaki1972} (in particular Theorem 3.1, Lemma 4.2, Theorem 4.1 and equation 6.11 in that paper) we will characterize $M_{t}^{Y}$ explicitly using the next set of arguments.

\begin{Def} Innovation process is defined as:
\begin{equation}
	d\nu^{Y}(t)=  dY(t)-E[h(X(t))\vert {\mathcal F}_{t}^{Y}]dt, \hspace{.1in} t\in [0,T]. \label{eqn5.7}
\end{equation}	
	
\end{Def} 
\begin{Lem} \cite{Fujisaki1972}  
	The innovation process $(\nu^{Y}, {\mathcal F}^{Y}_{t},m)$ is an $N$-vector standard Wiener process. Furthermore ${\mathcal F}^{Y}_{s}$ and $\sigma\left \{\nu_{v}^{Y}-\nu_{u}^{Y}; s\leq u\leq v \leq T\right \}$  are independent.
\end{Lem}

Let us now state an important martingale representation theorem (Theorem 3.1 \cite{Fujisaki1972}):
\begin{Thm}
	Every square integrable martingale $(M^{Y}_{t},{\mathcal F}^{Y}_{t},m)$ is sample continuous and has the representation
\begin{equation}
	M_{t}^{Y}-M^{Y}_{0}=\int_{0}^{t}\Xi(s)\cdot d\nu^{Y}(s).	\label{eqn5.8}
\end{equation}
where $\nu(t)$ is the innovation process and $\Xi(s)=(\Xi_{s}^1, \cdots,\Xi_{s}^{N})	$ is jointly measurable and adapted to ${\mathcal F}^{Y}_{s}$.	
\end{Thm}

\begin{Lem} \cite{Fujisaki1972}  Let $(M_{t}(f), {\mathcal F}_{t},m)$ be the square integrable martingale \ref{eqn5.5}. Then there exists unique sample continuous process $\langle M(f), Z^{i}\rangle, (i=1,\cdots,N)$ adapted to ${\mathcal F}_{t}$ such that almost all sample functions are of bounded variation and $M_{t}(f)-\langle M(f), Z^{i}\rangle_{t}$ are ${\mathcal F}_{s}$-martingales. Furthermore each $\langle M(f), Z^{i}\rangle_{t}$ has the following properties: it is absolutely continuous with respect to Lebesgue measure in $[0,T]$. There exists a modification of the Radon-Nikodym derivative which is $(t,\omega)$-measurable and adapted to ${\mathcal F}_{t}$ and which shall denote by $\tilde{D}_{t}^{i}f(\omega)$. Then using the vector notation $\tilde{D}_{t}f=(\tilde{D}_{t}^{1}f,\cdots, \tilde{D}_{t}^{N}f)$,
	\begin{displaymath}
		\langle M(f),Z \rangle_{t} =\int_{0}^{t}\tilde{D}_{s}fds, \hspace{.1in} \mbox{ a.s},
	\end{displaymath}
where
\begin{displaymath}
\int_{0}^{T}E\vert \tilde{D}_{s}f\vert^{2}ds <\infty.
\end{displaymath}
If the process $X(t)$ and $Z(t)$ are independent then 
\begin{displaymath}
	\langle M(f),Z \rangle_{t}=0, \mbox{ a.s}.
\end{displaymath}
	\end{Lem}

The following theorem characterizes $\Xi$ explicitly:
\begin{Thm}
Let $f\in D({\mathcal A})$ and 
\begin{equation}
	\int_{0}^{T}\vert h(X(t))\vert^{2}dt<\infty, \hspace{.1in} \mbox{ a.s.}\label{eqn5.10}
\end{equation}
Then the evolution of the conditional expectation $	E[f(X(t))\vert {\mathcal F}_{t}^{Y}]$ (the nonlinear filter) is characterized by the Fujisaki-Kallianpur-Kunita equation:
\begin{displaymath}
	E[f(X(t))\vert {\mathcal F}_{t}^{Y}]=	E[f(X(0))\vert {\mathcal F}_{0}^{Y}] + \int_{0}^{t}E[{\mathcal A}f(X(s))\vert {\mathcal F}_{s}^{Y}]ds
\end{displaymath}
\begin{equation}
+ \int_{0}^{t} \left \{E[f(X(s))h(X(s))\vert {\mathcal F}_{s}^{Y}] -E[f(X(s))\vert {\mathcal F}_{s}^{Y}]E[h(X(s))\vert {\mathcal F}_{s}^{Y}]	
	+ E[\tilde{D}_{s}f(X(s))\vert {\mathcal F}_{s}^{Y}]\right \}\cdot d\nu^{Y}(s) \label{eqn5.11}
\end{equation}
where  $\tilde{D}f(X(s))$ is given by
\begin{displaymath}
	\langle M(f),Z\rangle = \int_{0}^{t}\tilde{D}_{s}f ds.
\end{displaymath}
In particular $\tilde{D}_{t}f =0$ if $W$ and $Z$ are independent. Moreover, if 
\begin{displaymath}
	BW=\sigma_{1}W_{1}+\sigma_{2}Z,
\end{displaymath}
then 
\begin{displaymath}
	\tilde{D}_{t}f=\sigma_{2}^{*}\frac{\partial f}{\partial x}.
\end{displaymath}

\end{Thm}

Note that if $h$ has a growth condition
\begin{displaymath}
	\vert h(x)\vert \leq C \vert | x\vert |^{p}, \mbox{ for some } p\geq 1,  \hspace{.1in} x\in H,
\end{displaymath}
then by theorem \ref{Thm2.1} the condition \ref{eqn5.10} is satisfied.

\subsection{Existence and Uniqueness of Measure-valued solutions to Filtering Equations}

 A theorem of Getoor \cite{Getoor1975} provides the existence of a random measure $\pi_{t}^{Y}$ that is measurable with respect to ${\mathcal F}_{t}^{Y}$ such that 
 	\begin{displaymath}
 		E[f(X(t))\vert {\mathcal F}_{t}^{Y}]=\pi_{t}^{Y}[f] =\int_{H}f(\zeta)\pi_{t}^{Y}(d\zeta).
 	\end{displaymath}
   
 Substituting in \ref{eqn5.11}, the probability measure-valued process $\pi_{t}^{Y}[\cdot]$ satisfies the  Fujisaki-Kallianpur-Kunita (FKK) type equation: 
 	\begin{equation}
 		d\pi_{t}^{Y}[f]=\pi_{t}^{Y}[{\mathcal A}f]dt +\left(\pi_{t}^{Y}[hf+\tilde{D}_{t}f]-\pi_{t}^{Y}[h]\pi_{t}^{Y}[f]\right)\cdot d\nu^{Y}(t), \hspace{.1in} \mbox{ for } f\in {\mathcal D}_{\mbox{cyl}}. \label{eqn5.12}
 	\end{equation}	
  If we set 
 	\begin{equation}
 		\vartheta_{t}^{Y}[f]:= \pi_{t}^{Y}[f]\exp\left \{\int_{0}^{t}\pi_{t}^{Y}[h]\cdot dY(s) -\frac{1}{2}\int_{0}^{t}\vert \pi_{t}^{Y}[h]\vert^{2}ds\right \},\label{eqn5.13}
 	\end{equation}	
 then the measure-valued process  $\vartheta_{t}^{y}[\cdot]$ satisfies the Zakai type equation:
 	\begin{equation}
 		d\vartheta_{t}^{Y}[f]=\vartheta_{t}^{Y}[{\mathcal A}f]dt +\vartheta_{t}^{Y}[hf+\tilde{D}_{t}f]\cdot dY(t), \mbox{ for } f\in {\mathcal D}_{\mbox{cyl}}. \label{eqn5.14}
 	\end{equation}	
 Such measure-valued evolutions were first derived by Kunita \cite{Kunita1971} in the context of compact space valued signal processes.
 
One method of proving uniqueness of measure-valued solutions to filtering equations is to start with the unique solution of the backward Kolmogorov equation similar to the papers in the case of nonlinear filtering of stochastic Navier-Stokes equation \cite{Sritharan1995} or for the case of stochastic reaction diffusion equation \cite{Hobbs1996}.  Here we use the solution $\Phi(t)$ of the backward Kolmogorov equation with initial data $\Phi$ and express the Zakai equation as:
\begin{equation} 
	\vartheta_{t}^{Y}[\Psi]=\vartheta_{0}^{Y}[\Phi(t)]+ \int_{0}^{t} \vartheta_{s}^{Y}[ (\frac{\partial}{\partial s} +{\mathcal A})\Phi(s)]ds 
	+ \int_{0}^{t} \left(\vartheta_{s}^{Y}[(h+\tilde{D}_{t})\Phi(s)]\right)\cdot d\nu_{s},
\end{equation}	
and utilize the fact the $\Phi$ satisfies $(\frac{\partial}{\partial s} +{\mathcal A})\Phi(s)=0$ in the above equation (eliminating the integral involving the generator ${\mathcal A})$ and then proceed to prove the uniqueness of the random measures $\vartheta_{t}^{Y}$ as well as  $\pi_{t}^{Y}$. 

We now specialize to the case of signal noise $W$ and observation noise $Z$ being independent. In this setting we can use the results of Kunita \cite{Kunita1971} and Szpirglas\cite{Szpirglas1978} to obtain equivalent nonlinear filtering equation that do not explicitly involve the generator ${\mathcal A}$ of the signal Markov process $X$.

Uniqueness theorem relies on the following lemma \cite{Szpirglas1978, Kallianpur1984a}:
\begin{Lem}
Let us define a subclass of ${\mathcal D}({\mathcal A})$ as:
\begin{displaymath}
{\mathcal D}_{2}({\mathcal A}):= \left \{ f\in {\mathcal D}({\mathcal A}): {\mathcal A}f \in {\mathcal D}({\mathcal A}) \right \}.	
\end{displaymath}
 Let $\mu_{1},\mu_{2}\in {\mathcal M}(H)$ be finite measures on $(H, {\mathcal B}(H))$ such that
\begin{equation}
\langle \mu_{1}, f\rangle =\langle \mu_{2}, f\rangle, 	\hspace{.2in} \forall f\in {\mathcal D}_{2}({\mathcal A}),
\end{equation}
where $\langle \mu, f\rangle = \int_{H} f(x)d \mu(x).$ Then $\mu_{1}=\mu_{2}\in {\mathcal M}(H)$.
\end{Lem}
The uniqueness of measure-valued solutions for the nonlinear filtering equations and equivalence of two different forms, one involving the formal generator and other using Feller semigroup of the system Markov process eliminating the term involving the generator started with Kunita\cite{Kunita1971} and Szpirglass \cite{Szpirglas1978} and the following series of theorems can be proven by the same methods as in those original works for signal state space compact, locally compact and general Hilbert space cases as in \cite{Kunita1971, Stettner1989,Bhatt2000, Budhiraja2003,Fernando2013}: 
	
\begin{Thm}
	The random probability measure valued process $\pi^{Y}_{t} \in {\mathcal P}(H)$ uniquely solves the evolution equation:
		\begin{equation}
		\pi_{t}^{Y}[f]=\pi_{0}^{Y}[f]+\int_{0}^{t}\pi_{s}^{Y}[{\mathcal A}f]ds + \int_{0}^{t}\left(\pi_{s}^{Y}[hf]-\pi_{s}^{Y}[h]\pi_{s}^{Y}[f]\right)\cdot d\nu^{Y}(s), \hspace{.1in} \mbox{ for } f\in D({\mathcal A}).
		\end{equation}
	Equivalently, $\pi^{Y}_{t}$ uniquely solves the evolution equation:
\begin{equation} 
		\pi_{t}^{Y}[f]=\pi_{0}[P_{t}f]+ \int_{0}^{t} \left(\pi_{s}^{Y}((P_{t-s}f) h)-\pi_{s}^{Y}(P_{t-s}f) \pi_{s}^{Y}(h)\right)\cdot d\nu^{Y}_{s},	\hspace{.1in} \mbox{ for } f\in {\mathcal B}_{b}(H).
\end{equation}

\end{Thm}

\begin{Thm}
	the random measure valued process $\vartheta^{Y}_{t} \in {\mathcal M}(H)$ uniquely solves the evolution equation:
	\begin{equation}
		\vartheta_{t}^{Y}[f]=\vartheta_{0}^{Y}[f]+\int_{0}^{t}\vartheta_{s}^{Y}[{\mathcal A}f]ds + \int_{0}^{t}\vartheta_{s}^{Y}[hf]\cdot dY(s), \hspace{.1in} \mbox{ for } f\in D({\mathcal A}).
	\end{equation}
	Equivalently, $\vartheta^{Y}_{t} \in {\mathcal M}(H)$ uniquely solves the evolution equation:
	\begin{equation} 
		\vartheta_{t}^{Y}[f]=\vartheta_{0}[P_{t}f]+ \int_{0}^{t} \vartheta_{s}^{Y}((P_{t-s}f) h)\cdot dY({s}),	\hspace{.1in} \mbox{ for } f\in {\mathcal B}_{b}(H).
	\end{equation}
	
\end{Thm}

\begin{Thm}
Let $\pi^{Y}_{t}$ be the unique solution of the FKK equation for arbitrary initial condition $\pi_{0}$. Furthermore,
\begin{enumerate}
	\item The solution $\pi^{Y}_{t}$ is $\sigma(Y(s)-Y(0); 0\leq s\leq t)\vee\sigma(\pi_{0})$ -measurable.
	\item Let $\pi_{t}^{Y(\nu)}$ and $\pi_{t}^{Y(\mu)}$ be solutions with initial conditions $\pi_{0}=\nu$ and $\pi_{0}=\mu$ respectively, where $\mu,\nu \in {\mathcal M}(H)$.  
	Then for every $t>0$ 
	\begin{equation}
		\lim_{\nu\rightarrow \mu}E\left [ \vert\pi_{t}^{Y(\nu)}(f)-\pi_{t}^{Y(\mu)} (f)\vert \right] =0, \hspace{.1in} f \in C_{b}(H).
	\end{equation}
\end{enumerate}
\end{Thm}

\begin{Thm}
The filtering process $(\pi^{Y}, {\mathcal F}_{t}^{Y}, P_{\mu})$, $\mu \in {\mathcal M}(H)$ are Markov process associated with the transition probabilities $\Pi_{t}(\nu, \Gamma)$ defined by 
\begin{displaymath}
	\Pi_{t}(\nu,\Gamma)=P(\pi_{t}^{Y}\in \Gamma : \hspace{.1in} \Gamma \in {\mathcal B} ({\mathcal M}(H)),
\end{displaymath}
where ${\mathcal B} ({\mathcal M}(H))$ is the Borel algebra generated by the open (or closed) sets in ${\mathcal M}(H)$.  

Furthermore, the transition probabilities $\Pi_{t}(\nu,\Gamma)$ define a Feller semigroup in $C_{b}({\mathcal M}(H))$, where $C_{b}({\mathcal M}(H))$ is the space of all real continuous functions over ${\mathcal M}(H)$.
\end{Thm}
Systematic study of time asymptotic analysis of nonlinear filter was initiated by H. Kunita \cite{Kunita1971} and further developed by several authors \cite{Mazz1988, Stettner1989, Kunita1991, Ocone1996,Bhatt1995, Bhatt1999, Bhatt2000, Budhiraja2003}. We will further develop this subject to be applicable to classical and quantum systems studied in this paper.

\begin{Def}
	A probability measure $\nu\in {\mathcal P}(H)$ on the Borel algebra ${\mathcal B}(H)$ is called the barycenter of the probability measure $\Phi\in {\mathcal P}({\mathcal P}(H))$ on the Borel algebra ${\mathcal B}({\mathcal P}(H))$ if and only if for every $\phi\in C_{b}(H)$ 
	\begin{equation}
		\nu(\phi)=\int_{{\mathcal P}(H)}\nu'(\phi)\Phi(d\nu').	
	\end{equation}
\end{Def}
The following theorem is inspired by the original paper of Kunita\cite{Kunita1971} which was in turn technically corrected by several authors \cite{Baxendale2004,Budhiraja2003, VanHan2012}. 
\begin{Thm}
	Let us assume that the following condition regarding the sigma fields is satisfied:
	\begin{equation}
		\bigcap_{t\leq 0}\left ({\mathcal F}_{-\infty,0}^{Y}\vee {\mathcal F}^{X}_{-\infty, t}\right )={\mathcal F}_{-\infty,0}^{Y} \vee {\mathcal F}^{X}_{-\infty, -\infty}, \hspace{.1in} m \mbox{ a.s}.
	\end{equation}
	The existence and uniqueness of invariant measures for the transition probabilites $P_{t}(x,A), A \in {\mathcal B}(H)$, associated with the classical (Theorem 2.3, Theorem 3.5, Theorem 3.7) and quantum spin systems (Theorem 4.2) imply that the transition probabilities $\Pi_{t}(\nu,\Gamma)$ associated with the filtering process for these classican and quantum spin systems have unique invariant measures. The invariant measure $\mu$ for $P_{t}(x,A)$ is the barrycenter of the invariant measure $\Phi$ of $\Pi_{t}(\nu,\Gamma)$:
	\begin{equation}
		\mu(f)=\int_{{\mathcal M}(H)} \nu(f)\Phi(d\nu), \hspace{.1in} \forall f\in C_{b}(H).	
	\end{equation}
\end{Thm}

\subsection{Evolution Equation for Error Covariance }
The seminal paper of Kalman and Bucy \cite{Kalman1961} also gives a characterization of the evolution of error covariance using a Riccati equation. This concept was generalized by Kunita\cite{Kunita1991} for nonlinear filters. Let us describe the time evolution of error covariance slightly generalizing the results of Kunita\cite{Kunita1991}. The error covariance of general moments $f(X(t))$ and $g(X(t))$ is given by
\begin{equation}
	\macP_{t}[f,g]=E\left [\left (f(X(t))-\pi_{t}(f)\right )\left(g(X(t))-\pi_{t}(g)\right )\right].
\end{equation}
\begin{Thm} Let $f,g\in D({\mathcal A})$. Then $\macP_{t}[f,g]$ is differentiable with respect to $t$ and the derivative satisfies
	\begin{displaymath}
		\frac{d}{dt} \macP_{t}[f,g]= \macP_{t}[f,{\mathcal A}g] + \macP_{t}[{\mathcal A}f,g]+\macQ_{t}[f,g] 
	\end{displaymath}
\begin{equation}
	-\sum_{i=1}^{N} E \left [ \pi_{t} \left \{(f-\pi_{t}(f))(h^{i}-\pi_{t}(h^{i}))+\tilde{D}f\right \} \pi_{t}\left \{  (g-\pi_{t}(g))(h^{i}-\pi_{t}(h^{i}))+\tilde{D}g\right\}\right ],
\end{equation}
where
\begin{equation}
\macQ_{t}[f,g] =E\left [\mbox{ tr} (B^{*}Q B D_{x}f \otimes D_{x}g)\right].	
\end{equation}
\end{Thm}

{\bf Proof:}\\

Note that we have:
\begin{displaymath}
{\mathcal A}(fg)=f{\mathcal A}g + g{\mathcal A}f +\mbox{ tr} (B^{*}Q B D_{x}f \otimes D_{x}g).
\end{displaymath}
We first write the FKK equation for $\pi_{t}(fg)$ as
\begin{displaymath}
	\pi_{t}(fg)=\pi_{0}(fg) +\int_{0}^{t} \pi_{s}(f {\mathcal A}g)ds +\int_{0}^{t}\pi_{s}(g{\mathcal A}f)ds +\int_{0}^{s}\pi_{s}( \mbox{ tr} (B^{*}Q B D_{x}f \otimes D_{x}g))ds +M^{\pi}_{t},
\end{displaymath}
where $M^{\pi}_{t}$ is a ${\mathcal F}^{Y}_{t}$-martingale with mean zero and $f,g\in D({\mathcal A})$. Similarly since $\pi_{t}(f)$ and $\pi_{t}(g)$ satisfy FKK equation we have by Ito formula for the product $\pi_{t}(f)\pi_{t}(g)$:
\begin{displaymath}
	\pi_{t}(f)\pi_{t}(g)=\pi_{0}(f)\pi_{0}(g)+ \int_{0}^{t}\pi_{s}(f)\pi_{s}({\mathcal A}g)ds+\int_{0}^{t}\pi_{s}(g)\pi_{s}({\mathcal A}f)ds
\end{displaymath}
\begin{displaymath}
+ \sum_{i}\int_{0}^{t} \left \{ \pi_{s}(fh^{i})-\pi_{s}(f)\pi_{s}(h^{i})+\pi_{s}(\tilde{D}f)\right \}\left \{ \pi_{s}(gh^{i})-\pi_{s}(g)\pi_{s}(h^{i})+\pi_{s}(\tilde{D}g)\right \}ds +\tilde{M}^{\pi}_{t},
\end{displaymath}
where $\tilde{M}^{\pi}_{t}$ is a ${\mathcal F}^{Y}_{t}$-martingale with mean zero. Now noting that
\begin{displaymath}
	\pi_{s}(fg)-\pi_{s}(f)\pi_{s}(g)=\pi_{s}\left \{(f-\pi_{s}(f))(g-\pi_{s}(g))\right\},
\end{displaymath}
we get 
\begin{displaymath}
	\macP_{t}(f,g)=E[\pi_{t}(fg)-\pi_{t}(f)\pi_{t}(g)].
\end{displaymath}
Substituting for $\pi_{t}(fg)$ and $\pi_{t}(f)\pi_{t}(g)$ from  the previous two equations and taking expectation we arrive at
\begin{displaymath}
	\macP_{t}(f,g)=\macP_{0}(f,g) +\int_{0}^{t}\left \{ \macP_{s}[f,{\mathcal A}g] + \macP_{s}[{\mathcal A}f,g]+\macQ_{s}[f,g] \right \}ds
\end{displaymath}
\begin{displaymath}
	-\sum_{i=1}^{N}\int_{0}^{t} E \left [ \pi_{s} \left \{(f-\pi_{s}(f))(h^{i}-\pi_{s}(h^{i}))+\tilde{D}f\right \} \pi_{s}\left \{  (g-\pi_{s}(g))(h^{i}-\pi_{s}(h^{i}))+\tilde{D}g\right\}\right ]ds.	
\end{displaymath}

 \section{White Noise Theory of Nonlinear Filtering}
 White noise theory of nonlinear filtering was initiated by several papers of A. V. Balakrishnan \cite{Bala1980, Bala1976, Bala1980} and systematically developed by G. Kallianpur and R. Karandikar \cite{Kallianpur1983, Kallianpur1984a, Kallianpur1984b, Kallianpur1985, Kallianpur1988} (see also Kallianpur's review \cite{Kallianpur1991}). We will further develop this method to be applicable to classical and quantum spin systems discussed in this paper. Key mathematical equations of white noise nonlinear filtering theory are essentially parallel to stochastic calculus counterpart but of deterministic in nature. The Kallianpur-Striebel formula, the Fujisaki-Kallianpur-Kunita equation, the Zakai equation and the Kunita's semigroup versions of FKK and the Zakai equation all have white noise counterparts and will be developed below for the spin system models. We will also establish the equivalency of FKK, Zakai with their Kunita counter parts in the spirit of Szpirglas\cite{Szpirglas1978}. 
 
 We will now describe the theory of Segal \cite{Segal1956} and Gross \cite{Gross1960, Gross1962} of finitely additive cylindrical measures on separable Hilbert spaces. Let ${\mathcal H}$ be a separable Hilbert space and ${\mathcal P}$ be the set of orthogonal projections on ${\mathcal H}$ having finite dimensional range. For $P\in {\mathcal P}$, let ${\mathcal C}_{P}=\{P^{-1}B: B \mbox{ a Borel set in the range of } P\}$. Let ${\mathcal C}=\cup_{P\in {\mathcal P}} {\mathcal C}_{P}$. A cylinder measure $\n$ on ${\mathcal H}$ is a finitely additive measure on $({\mathcal H}, {\mathcal C})$ such that its restriction to ${\mathcal C}_{P}$ is countably additive for each $P\in {\mathcal P}$.
 
 Let $L$ be a representative of the weak-distribution corresponding to the cylinder measure $\n$. This means that $L$ is a linear map from ${\mathcal H^{*}}$ (identified with ${\mathcal H}$) in to ${\mathcal L}_{0}(\Omega_{1}, \Sigma_{1},m_{1})$ -the space of all random variables on a $\sigma$-additive probability space $(\Omega_{1}, \Sigma_{1},m_{1})$ such that 
 \begin{displaymath}
 	\n(h\in {\mathcal H}: ((h,h_{1}), (h,h_{2}),\cdots, (h,h_{k}))\in B)
\end{displaymath}
\begin{equation}
 	= m_{1}(\omega\in \Omega_{1}; \left(L(h_{1})(\omega), L(h_{2})(\omega), \cdots, L(h_{k})(\omega)\right)\in B) \label{eqn6.1}
 \end{equation}
for all Borel sets $B$ in $\mR^{k}$, $h_{1}, \cdots, h_{k}\in {\mathcal H}$, $k\geq 1$. Two maps $L, L'$ are said to be equivalent if both satisfy \ref{eqn6.1} and the equivalence class of such maps is the weak distribution corresponding to the cylindrical measure $\n$.
 
 A function $f$ on ${\mathcal H}$ is called a tame function if it is of the form
 \begin{equation}
 	f(y)=\varphi((y,h_{1}),\cdots, (y,h_{k})), \mbox{ for some } k\geq 1, h_{1}, \cdots, h_{k}\in {\mathcal H}  \label{eqn6.2}
 \end{equation}	
and a Borel function $\varphi:\mR^{k}\rightarrow \mR$. For a tame function $f$ given by \ref{eqn6.2}, we associate the random variable $\varphi(L(h_{1}), \cdots, L(h_{k}))$ on $(\Omega_{1},\Sigma_{1},m_{1})$ and denote it by $\tilde{f}$. We can extend this map $f\rightarrow \tilde{f}$ to a larger class of functions as follows \cite{Gross1962}:
\begin{Def} \label{D6.1} Let ${\mathcal L}({\mathcal H}, {\mathcal C}, \n)$ be the class of continuous functions $f$ on ${\mathcal H}$ such that the net $\{\tilde{(f\circ P)}: P\in {\mathcal P}\}$ is Cauchy in $m_{1}$-measure. Here $P_{1}<P_{2}$ if the range $P_{1} \subseteq$ range $P_{2}$. For $f\in {\mathcal L}({\mathcal H}, {\mathcal C}, \n)$, we define:
	\begin{displaymath}
	\tilde{f}={\mbox{Lim in Prob}}_{P\in {\mathcal P}} \tilde{(f\circ P)	}.
	\end{displaymath}
The map $f\rightarrow \tilde{f}$ is linear, multiplicative and depends only on $f$ and $\n$ and is independent of the representative of the weak distribution.	
\end{Def}

\begin{Def} \label{D6.2} The function $f\in {\mathcal L}({\mathcal H}, {\mathcal C}, \n)$ is integrable with respect to $\n$ if $\int_{\Omega_{1}}\vert \tilde{f}\vert dm_{1}<\infty$ and then for $C\in {\mathcal C}$ define the integral $f$ with respect to $\n$ over $C$ denoted by $\int_{C}fd\n$ by
	\begin{displaymath}
\int_{C}fd\n =\int_{\Omega_{1}} \tilde{1_{C}} \tilde{f} dm_{1}.		
	\end{displaymath}
\end{Def}	
Let us now consider the class of all Gauss measures on ${\mathcal H}$ defined by
\begin{displaymath}
	\mu\{y\in {\mathcal H}: (y,h)\leq a\}=\frac{1}{\sqrt{2\pi v(h)}}\int_{-\infty}^{a} \exp{(-\frac{x^{2}}{2v(h)})}dx, \hspace{.1in} \forall h\in {\mathcal H}.
\end{displaymath}
The special case of $v(h)=\|h\|^{2}$ is called the canonical Gauss measure $\m$. The following Lemma from Sato \cite{Sato1969} (in particular, Lemma 6 in that paper)clarifies the $\sigma$-additivity of Gauss measures on separable Hilbert spaces.	
\begin{Lem} \label{L6.1} let ${\mathcal H}$ be a separable Hilbert space and let $\mu$ be a Gaussian cylinder measure on $({\mathcal H}, {\mathcal C})$ with variance $v(h)$. Then $\mu$ has a $\sigma$-additive extension to $({\mathcal H}, \bar{\mathcal C})$ (here $\bar{\mathcal C}$ is the minimal $sigma$-algebra containing ${\mathcal C})$ if and only if the characteristic functional of $\mu$ is of the form:
	\begin{displaymath}
	\int_{\mathcal H}e^{i (h,x)}d\mu(x)	=\exp {\left [ i<h,m_{e}>-\frac{1}{2}\|Sh\|^{2}\right]}, \hspace{.1in} \forall h\in {\mathcal H},
	\end{displaymath}
where $	m_{e}$ is an element of ${\mathcal H}$, and $S$ is a nonnegative selfadjoint Hilbert-Schmidt operator on ${\mathcal H}$.
\end{Lem}
The canonical Gauss measure $\m$ corresponds to $m_{e}=0$ and $S=I$ and since the identity operator $I$ is not Hilbert-Schmidt, the cannonical Gauss measure $\m$ is only finitely additive (See also Chapter-I, Proposition 4.1 of Kuo\cite{Kuo1975}). 

The identity map $e$ on ${\mathcal H}$, considered as a map from $({\mathcal H}, {\mathcal C}, \n)$ in to $({\mathcal H}, {\mathcal C})$ is called the Gaussian white noise.

Let us now start with the abstract version of the white noise filtering model as formulated by Kallianpur and Karandikar \cite{Kallianpur1988}:
\begin{equation}
	y=\zeta + e \label{eqn6.3}
\end{equation}
where $\zeta$ is an ${\mathcal H}$-valued random variable defined on a countably additive probability space $(\Omega, {\mathcal F}, m)$, independent of $e$. To formulate this mathematicially precisely let ${\mathcal E}={\mathcal H}\times \Omega$ and
\begin{displaymath}
	\Sigma=\cup_{P\in {\mathcal P}} {\mathcal C}_{P}\otimes {\mathcal F}  \mbox{ where } {\mathcal C}_{P}\otimes {\mathcal F} \mbox{ is the product sigma field}.
\end{displaymath}
For $P\in {\mathcal P}$, let $\alpha_{P}$ be the product measure $(\m\vert {\mathcal C}_{P})\otimes m$ which is countably additive. This defines a unique finitely additive probability measure $\alpha$ on $({\mathcal E}, \Sigma)$ such that $\alpha=\alpha_{P}$ on ${\mathcal C}_{P}\otimes {\mathcal F} $.

let $e, \zeta, y$ be $H$-valued maps on ${\mathcal E}$ defined by
\begin{displaymath}
	e(h, \omega)=h,
	\end{displaymath}
\begin{displaymath}
	\zeta(h,\omega)=\zeta(\omega),
\end{displaymath}
\begin{equation}
y(h,\omega)=e(h,\omega)+ \zeta(h,\omega), \hspace{.1in} (h,\omega)\in {\mathcal H}\times \Omega. \label{eqn6.4}
\end{equation}
The model \ref{eqn6.4} is the abstract version of the white noise filtering on $({\mathcal E}, \Sigma, \alpha)$. Our goal is to first characterize the conditional expectation $E[g\vert y]$ in this finitely additive setting. 
\begin{Def}
If there exists a $v\in {\mathcal L}({\mathcal H}, {\mathcal C}, \n)$ such that 
\begin{equation}
\int_{{\mathcal H}\times \Omega}g(\omega) 1_{C}(y(h,\omega))d\alpha(h,\omega)=\int_{C}v(y)d\n(y),\hspace{.1in} \forall C\in {\mathcal C},\label{eqn6.5}
\end{equation}then we define $v$ to be the conditional expectation of $g$ given $y$ and express it as $E[g\vert y]=v.$
\end{Def}
Note that the integrand $g(\omega) 1_{C}(y(h,\omega))$ is ${\mathcal C}_{P}\otimes {\mathcal F}$-measurable for $C\in {\mathcal C}_{P}$ and hence the integrals in \ref{eqn6.5} are well-defined.
Let us now state the elegant Bayes formula for the finitely additive white noise framework proved in \cite{Kallianpur1983}:
\begin{Thm} 
	Let $y, \zeta$ be as in \ref{eqn6.4}. Let $g$ be an integrable function on $(\Omega, {\mathcal F}, m)$. Then 
\begin{equation}
E[g\vert y] =\frac{\int_{\Omega}g(\omega)exp{\left \{(y,\zeta(\omega))-\frac{1}{2}\|\zeta(\omega)\|^{2}\right\}}dm(\omega)}{\int_{\Omega}exp{\left \{(y,\zeta(\omega))-\frac{1}{2}\|\zeta(\omega)\|^{2}\right\}}dm(\omega)}. \label{eqn6.6}
\end{equation}
\end{Thm}
We now return to our spin system models where the state variable $X(t)$ is an $(H, {\mathcal B}(H))$-valued stochastic process (where ${\mathcal B}(H)$ is the Borel algebra generated by the closed or open subsets of $H$) that is defined on the probability space $(\Omega, {\mathcal F},m)$. Let ${\mathcal K}$ be a separable Hilbert space (including $\mR^{N}$) and let the observation vector $\zeta=\zeta(X(t))$ be a measurable function from $H\rightarrow {\mathcal K}$ such that 
\begin{equation}
	\int_{0}^{T}\|\zeta(X(t))\|_{\mathcal K}^{2}dt<\infty, \mbox{ a.s.}.
\end{equation}
The Hilbert space ${\mathcal H} $ we started with in the abstract white noise formulation will be in this case ${\mathcal H}=L^{2}(0,T;{\mathcal K})$.

 The white noise sensor measurement model:
 	
 	\begin{equation}
 		Y(t) =\zeta(X(t))+ e(t) \in {\mathcal K},\hspace{.1in} t>0,
 	\end{equation}
 	where $e(t)\in {\mathcal K}$ is a  finite or infinite dimensional white noise. 
 	
 	Let us define a Borel measure $\rho^{Y}_{t}(\cdot)\in {\mathcal M}(H)$ and a probability measure $\pi_{t}^{Y}(\cdot)\in {\mathcal P}(H)$ as follows:
 	\begin{displaymath}
 	\rho^{Y}_{t}(B)	=E\left [ 1_{B}(X(t))\exp{\int_{0}^{t}\left \{(\zeta(X(s)), Y(s))_{\mathcal K}-\frac{1}{2}\|\zeta(X(s))\|^{2}_{\mathcal K} \right \}ds}\right ],
 	\end{displaymath}
for $B\in {\mathcal B}(H)$, and
 \begin{displaymath}
 	\pi^{Y}_{t}(B)=\frac{\rho^{Y}_{t}(B)}{\rho^{Y}_{t}(H)}.
 \end{displaymath}
 	Then the measures $\rho^{Y}_{t}\in {\mathcal M}(H)$ and $\pi_{t}^{Y}\in {\mathcal P}(H)$ satisfy: 
 	\begin{equation}
 		<\pi_{t}^{Y},f>=\int_{H}f(x)\pi^{Y}_{t}(dz)=E\left [ f(X(t))\vert Y(s), 0\leq s\leq t \right ],
 	\end{equation}
 
 	\begin{displaymath}
 		<\pi_{t}^{Y}, f> = \frac{<\rho^{Y}_{t}, f>}{<\rho_{t}^{Y},1>},
 	\end{displaymath}
 
 	\begin{displaymath}
 		<\rho_{t}^{Y}, f>= E\left \{ f(X(t)) \exp{\int_{0}^{t}{\mathfrak{C}}_{s}^{Y}(X(s))ds }\right \},
 	\end{displaymath}
 
 where  
 	\begin{displaymath}
 		{\mathfrak C}_{s}^{Y}(X)=(\zeta(X(s)), Y(s))_{\mathcal K}-\frac{1}{2}\|\zeta(X(s))\|^{2}_{\mathcal K}.
 	\end{displaymath}
 
The following series of theorems on equivalence of finitely additive nonlinear filtering equations, uniqueness of measure-valued filter, Markov property and robustness are proven by methods similar to those of Kallianpur and Karandikar \cite{Kallianpur1983, Kallianpur1984a, Kallianpur1984b, Kallianpur1988, Kallianpur1991}:

\begin{Thm}
	For $Y\in C([0,T];{\mathcal H})$ and the class of measures $\rho_{t}^{Y}\in {\mathcal M}(H)$ uniquely solves
		\begin{equation}
			<\rho_{t}^{Y},f>=<\rho_{0},f>+\int_{0}^{t}<\rho_{s}^{Y},{\mathcal A}f+{\mathfrak C}_{s}^{Y}f>ds, \hspace{.1in} f\in D({\mathcal A}),0\leq t\leq T,
		\end{equation}	
or equivalently, $\rho_{t}^{Y}\in {\mathcal M}(H)$ uniquely solves	
	\begin{equation}
	<\rho_{t}^{Y},f>=<\rho_{0},P_{t}f>+\int_{0}^{t}<\rho_{s}^{Y},{\mathfrak C}_{s}^{Y}(P_{t-s}f)>ds, \hspace{.1in} f\in {\mathcal B}_{b}(H), 0\leq t\leq T.
\end{equation}	

\end{Thm}

\begin{Thm}
	For $Y\in C([0,T];{\mathcal H})$ the probability measure valued process $\pi_{t}^{Y}\in {\mathcal P}(H)$, $0\leq t\leq T$, uniquely solves
	\begin{equation}
		<\pi_{t}^{Y},f>=<\pi_{0},f>	+\int_{0}^{t}\left [<\pi_{s}^{Y},{\mathcal A}f 	+{\mathfrak C}_{s}^{Y}f>-<\pi_{s}^{Y},{\mathfrak C}_{s}^{Y}><\pi_{s}^{Y},f>\right ] ds,  f\in D({\mathcal A}),
	\end{equation}
or equivalently $\pi_{t}^{Y}\in {\mathcal P}(H)$ uniquely solves	
	\begin{equation}
	<\pi_{t}^{Y},f>=<\pi_{0},P_{t}f>		+\int_{0}^{t}\left [<\pi_{s}^{Y},{\mathfrak C}_{s}^{Y}(P_{t-s}f)>-<\pi_{s}^{Y},{\mathfrak C}_{s}^{Y}><\pi_{s}^{Y},(P_{t-s}f)>\right ] ds, f\in {\mathcal B}_{b}(H).
	\end{equation}
\end{Thm}	

Let
\begin{displaymath}
	{\mathcal H}_{t}=\left \{ \eta \in {\mathcal H}: \int_{t}^{\infty}\|\eta(r)\|^{2}_{\mathcal K}dr=0\right \}.
\end{displaymath}
Note that ${\mathcal H}_{t}$ is a closed subspace of ${\mathcal H}$. Let $Q_{t}$ be the orthogonal projection onto ${\mathcal H}_{t}$ and let ${\mathcal C}_{Q_{t}}={\mathcal C}({\mathcal H}_{t})$.

\begin{Thm}
$\pi^{Y}_{t}$ and $\rho^{Y}_{t}$ are respectively ${\mathcal P}(H)$ and ${\mathcal M}(H)$ -valued $\{{\mathcal C}_{Q_{t}}\}$ Markov processes on $({\mathcal H}, {\mathcal C}, \n)$.
\end{Thm}

\begin{Thm}
If $Y_{n}\rightarrow Y$ in ${\mathcal H}$ then $\pi_{t}^{Y_{n}}\rightarrow \pi_{t}^{Y}$ and $\rho_{t}^{Y_{n}}\rightarrow \rho_{t}^{Y}$ in total variation norm.
\end{Thm}

\vspace{.3in}

{\bf Acknowledgment}: The first author's research has been supported by the U. S. Air Force Research Laboratory through the National Research Council Senior Research Fellowship of the National Academies of Science, Engineering and Medicine.



\begin{thebibliography}{30} 
\bibitem{Albeverio2001} S. Albeverio, Yu G. Kondratiev, M. Rockner and T. V. Tsikalenko, Gauber dynamics of quantum lattice systems, Reviews in Mathematical Physics, Vol. 13, No. 1 (2001) 51-124.	
\bibitem{Arous1995} G. B. Arous and A. Guionnet, Large deviations for Langevin spin glass dynamics, {\it Probab. Theory Relat. Fields} Vol. 102, (1995), 455- 509.
\bibitem{Bala1974} A. V. Balakrishnan, On the approximation of It\^{o} integrals using band-limited processes. {\it SIAM J. Control}, 12:237–25, 975, 1974.
\bibitem{Bala1976}  A. V. Balakrishnan, Radon-Nikodym derivatives of a class of weak distributions on Hilbert spaces, {\it Appl. Math. Optim.}, 3(2-3):209–225, 1976/77.
\bibitem{Bala1980} A. V. Balakrishnan, Nonlinear white noise theory. In {\it Multivariate analysis}, V (Proc. Fifth Internat. Sympos., Univ. Pittsburgh, Pittsburgh, Pa., 1978), pages 97–109. North-Holland, Amsterdam, (1980).
\bibitem{Bhatt1995} A. G. Bhatt, G. Kallianpur, and R. Karandikar, Uniqueness and robustness of solution of measure-valued equations of nonlinear filtering, {The Annals of Probability}, Vol. 23, No.4, (1995),1895-1938.
\bibitem{Baxendale2004} P. Baxendale, P. Chigansky, and R. Liptser Asymptotic Stability of the Wonham Filter: Ergodic and Nonergodic Signals, {\it SIAM Journal of Control and Optimization}, Vol. 4, Issue 2, (2004), 643-669.
\bibitem{Belavkin1980} V. P. Belavkin, Quantum filtering of Markov signals with white quantum noise, {\it Radiotechnika i Electronika}, 25 (1980), pp. 1445–1453.
\bibitem{Belavkin1983} V. P. Belavkin, Theory of the control of observable quantum systems, {\it Autom. Rem. Control,}44 (1983), pp. 178–188.
\bibitem{Belavkin1987} V. P. Belavkin, Nondemolition measurement and control in quantum dynamical systems, in {\it Information Complexity and Control in Quantum Physics}, CISM Courses and Lectures 294, S. Diner and G. Lochak, eds., Springer-Verlag, Vienna, 1987, pp. 331–336.
\bibitem{Belavkin1992a} V. P. Belavkin, Quantum continual measurements and a posteriori collapse on CCR, {\it Comm. Math. Phys.}, 146 (1992), pp. 611–635.
\bibitem{Belavkin1992b} V. P. Belavkin, Quantum stochastic calculus and quantum nonlinear filtering, {\it J. Multivariate Anal.}, 42 (1992), pp. 171–201.
\bibitem{Bhatt1999} A. G. Bhatt, G. Kallianpur, and R. Karandikar, Robustness of nonlinear filter, {Stochastic Processes and their Applications}, 81 (1999) 247-254.
\bibitem{Bhatt2000} A. G. Bhatt, A. Budiraja, and R. Karandikar, Markov property and ergodicity of the nonlinear filter, {SIAM J. of Control and Optimization}, Vol. 39, No. 3, (2000), pp. 928–949.
\bibitem{Bouten2007} L. Bouten, R. van Handel and M. R. James, An introduction to quantum filtering, {\it SIAM J. Control and Optimization}, Vol. 46, No. 6, (2007), pp. 2199-2241.
\bibitem{Bratteli1987} O. Bratteli and D. W. Robinson, {\it Operator Algebras and Quantum Statistical Mechanics-I}, Springer-Verlag, Berlin (1987).
\bibitem{Bratteli1996} O. Bratteli and D. W. Robinson, {\it Operator Algebras and Quantum Statistical Mechanics-II}, Springer-Verlag, Berlin (1996).
\bibitem{Budhiraja2003} A. Budhiraja, Asymptotic stability, ergodicity and other asymptotic properties of the nonlinear filter, {\it Ann. I. H. Poincaré – PR}, 39, 6 (2003) 919–941.
\bibitem{DaPrato1995} G. Da Prato and J. Zabczyk, Convergence to equilibrium for classical and quantum spin systems, {\it Probab. Theory Relat. Fields}, Vol. 103, (1995), 529-552.
\bibitem{DaPrato1996} G. Da Prato and J. Zabczyk, {\it Ergodicity for Infinite Dimensional Systems}, London Mathematical Society Lecture Note Series. 229, (1996).
\bibitem{DaPrato2014} G. Da Prato and J. Zabczyk, {\it Stochastic Equations in Infinite Dimensions},Second Edition,  Cambridge University Press, (2014).
\bibitem{Fernando2013} B. P. W. Fernando and S. S. Sritharan, Nonlinear Filtering of Stochastic Navier-Stokes equation with Ito-Levy noise, {\it Stochastic Analysis and Applications}, Vol.  31, No. 3, (2013), 381-426.
\bibitem{Fujisaki1972} M. Fujisaki, G. Kallianpur, and H. Kunita, Stochastic differntial equations for the nonlinear filtering problem, {\it Osaka J. Math.} (1972), Vol.9, 19–40.
\bibitem{Getoor1975}  R. K. Getoor,  On the construction of kernels. In Seminaire de probabilitiés IX. {\it Lecture Notes in Mathematics.} Meyer, P.A., ed., Vol. 465, Springer-Verlag, Berlin (1975), 443-463.
\bibitem{Gross1960} L. Gross, Integration and nonlinear transformation in Hilbert space, {\it Trans. Amer. Math. Soc. }, Vol. 94, (1960), 404-440.
\bibitem{Gross1962} L. Gross, Measurable functions on Hilbert space, {\it Trans. Amer. Math. Soc.}, Vol. 105, (1960), 372-390.
\bibitem{Hobbs1996} S. Hobbs and S. S. Sritharan, Nonlinear filtering for stochastic reaction and diffusion equations, in {\it Probability Theory and Modern Analysis},  N. Gretsky, J. Goldstein and J. J. Uhl, eds, Marcel Dekker, New York, (1996).
\bibitem{Holley1981}R. Holley and D. Stroock, Diffusions on an infinite dimensional torus, {\it J. of. Functional Analysis}, Vol. 42, (1981), 29-63.
\bibitem{Hudson1984} R. L. Hudson and K. R. Parthasarathy, Quantum Ito's Formula and Stochastic Evolutions.{\it Comm. Math. Phys.} 93, 301-323, (1984).
\bibitem{Kallianpur1969} G. Kallianpur and C. Striebel, Stochastic differential equations occuring in the estimation of continuous parameter stochastic processes, {\it The Probability Theory and its Applications}, Vol. XIV, No.4, (1969), 567-594.
\bibitem{Kallianpur1983} G. Kallianpur and R. Karandikar, A Finitely Additive White Noise Approach to Nonlinear Filtering, {\it Appl. Math. Optim.} Vol.10, (1983), 159-185.
\bibitem{Kallianpur1984a} G. Kallianpur and R. Karandikar, Measure-Valued Equations for the Optimum Filter in Finitely Additive Nonlinear Filtering Theory, {\it Z. Wahrscheinlichkeitstheorie verw. Gebiete}, 	Vol. 66, (1984), 1-17.
\bibitem{Kallianpur1984b} G. Kallianpur and R. Karandikar, The Markov Property of the Filter in the Finitely Additive White Noise Approach to Nonlinear Filtering, {\it Stochastics}, Vol. 13, (1984), 177-198.
\bibitem{Kallianpur1985} G. Kallianpur and R. Karandikar, White noise calculus and nonlinear filtering theory, {\it The Annals of Probability}, Vol. 13, No. 4, (1985), 1033-1107.
\bibitem{Kallianpur1988} G. Kallianpur and R. Karandikar, {\it White Noise Theory of Prediction, Filtering and Smoothing}, Gordon and Breach Publishers, Amsterdam (1988).
\bibitem{Kallianpur1991} G. Kallianpur, A skeletal theory of filtering, {Stochastic Analysis}, Edited by E. Mayer-Wolf, E. Merzbach, and A. Shwartz, Academic Press (1991).
\bibitem{Kalman1961} R. E. Kalman and R. S. Bucy, New results in linear filtering and prediction theory, {Journal of Basic Engineering}, March (1961), 95-108.
\bibitem{Kunita1971} H. Kunita, Asymptotic behavior of nonlinear filtering error of Markov processes, {\it J. of Multivariate Analysis}, vol.1, (1971), 365-393.
\bibitem{Kunita1991} H. Kunita, Ergodic properties of nonlinear filtering processes, in K. S. Alexander and J. Watkins editors {\it  Spatial stochastic processes : a festschrift in honor of Ted Harris on his seventieth birthday }, Springer-Verlag (1991), 234-259.
\bibitem{Kuo1975} H-H. Kuo, {\it Gaussian Measures in Banach Spaces}, Springer-Verlag, New York, (1975).
\bibitem{Mazz1988} G. Mazziotto, L. Stettner, J. Szpirglas and J. Zabczyk, On impulse control with partial observation, {SIAM J. of Control and Optimization}, Vol. 26, No.4, (1988), 964-984.
\bibitem{Ocone1996} D. Ocone and E. Pardoux, Asymptotic stability of the optimal filter with respect to its initial condition, {\it SIAM J. of Control and Optimization}, Vol. 34, No. 1, (1996), 226-243.
\bibitem{Partha1992} K. R. Parthasarathy, {\it An Introduction to Quantum Stochastic Calculus}, Birkauser, Boston, (1992).
\bibitem{Peszat2007} S. Peszat and J. Zabczyk, {\it Stochastic Partial Differential Equations with L\'evy Noise}, Cambridge University Press, (2007).
\bibitem{Rebes2015} P. Rebeschini and R. van Handel,  Phase transitions in nonlinear filtering, {\it     Electron. J. Probab.} 20 (2015), no. 7, 1–46.
\bibitem{Ruelle1969} D. Ruelle, {\it Statistical Mechanics: Rigorous Results}, World Scientific., Singapore (1969).
\bibitem{Sato1969} H. Sato, Gaussian measure on a Banach space and abstract Wiener measure, {\it Nagoya Math. J.} Vol. 36, (1969), 65-81.
\bibitem{Segal1956} I. E. Segal, Tensor algebras over Hilbert spaces, {Trans. Amer. Math. Soc.} Vol. 81, (1956), 106-134.
\bibitem{Simon1974} B. Simon, The $P(\Phi)_{2}$-Euclidean (Quantum) Field Theory, Princeton University Press, Princeton (1974).
\bibitem{Simon1993} B. Simon, Statistical Mechanics of Lattice Gases, Princeton University Press, Princeton (1993).
\bibitem{Sritharan1995}  S. S. Sritharan, Nonlinear filtering of stochastic Navier-Stokes equations. In {\it Nonlinear Stochastic PDEs: Burgers Turbulance and Hydrodynamic Limit},T. Funaki, T., and W. A. Woyczynski eds., Springer-Verlag, New York, (1995) 247–260.
\bibitem{Sritharan2010} S. S. Sritharan and M. Xu, Convergence of particle filtering method for nonlinear estimation of vortex dynamics. {\it Comm. Stoch. Anal.} Vol. 4, No. 3, 443–465, (2010).
\bibitem{Sritharan2011} S. S. Sritharan and M. Xu, A stochastic Lagrangian particle model and nonlinear filtering for three dimensional Euler flow with jumps, {\it Comm. Stoch. Anal.}, Vol. 5, No.3, (2011), 565–583.
\bibitem{Stettner1989} L. Stettner, On invariant measures of filtering processes, {\it Lecture Notes in Control and Inform. Sci.}, Vol. 126, Springer, Berlin, 1989, 279–292.
\bibitem{Szpirglas1978} J. Szpirglas, Sur l’équivalence d’équations différentielles stochastiques à valeurs mesures intervenant dans le filtrage markovien non linéaire, {\it Ann. Inst. Henri Poincaré,} 	Vol. XIV, n° 1, (1978), 33-59. 
\bibitem{VanHan2012} R. van Handel, On the exchange of intersection and supremum of $\sigma$-fields in filtering theory, {\it Israel Journal of Mathematics} volume 192,  (2012), 763-782.
\bibitem{Zagarlinski1996} B. Zegarlinski, The strong decay to equilibrium for the stochastic dynamics of unbounded spin systems on a lattice, {\it Commun. Math. Phys.} 175, 401-432 (1996).
\bibitem{Zakai1969} M. Zakai, On the optimal filtering of diffusion processes, {\it Wahrsh, Z. Verw. Gebiete}, (1969), Vol. 11, 230–243.

\end{thebibliography}
\end{document}